\documentclass[12pt]{article}
\usepackage{amsmath, amsthm, amssymb, graphicx, array, booktabs, caption}
\usepackage{dsfont}
\usepackage{hyperref}
\usepackage{fullpage}

\hypersetup{
    colorlinks=true,
    linkcolor=blue,
    filecolor=magenta,
    citecolor=green,    
    urlcolor=cyan,
    pdfpagemode=FullScreen,
    }
    
\begin{document}

\newcommand{\Z}{\mathbb{Z}}
\newcommand{\N}{\mathbb{N}}
\newcommand{\ind}{\mathds{1}}
\newcommand{\Topp}[2]{\mathrm{Topp}^{#1}_{#2}}
\newcommand{\Stab}{\mathrm{Stab}}
\newcommand{\Det}{\mathrm{Deterministic}}
\newcommand{\A}{\mathrm{Always}\text{-}}
\newcommand{\FT}{\mathrm{FT}}
\newcommand{\NFT}{\mathrm{NFT}}

\hyphenation{
  neigh-bou-ring
  cu-mu-la-tive
  dis-tri-bu-tion
  dis-tri-bu-tions
  Equa-tion
  para-me-ters
}  

\title{Asymptotics of the single-source stochastic sandpile model}
\author{Thomas Selig and Haoyue Zhu}
\date{\today}

\maketitle

\abstract{In the single-source sandpile model, a number $N$ grains of sand are positioned at a central vertex on the 2-dimensional grid $\mathbb{Z}^2$. We study the stabilisation of this configuration for a stochastic sandpile model based on a parameter $M \in \mathbb{N}$. In this model, if a vertex has at least $4M$ grains of sand, it topples, sending $k$ grains of sand to each of its four neighbours, where $k$ is drawn according to some random distribution $\gamma$ with support $\{1,\cdots,M\}$. Topplings continue, a new random number $k$ being drawn each time, until we reach a stable configuration where all the vertices have less than $4M$ grains. This model is a slight variant on the one introduced by Kim and Wang.

We analyse the stabilisation process described above as $N$ tends to infinity (with $M$ fixed), for various probability distributions $\gamma$. We focus on two global parameters of the system, referred to as radius and avalanche numbers. The radius number is the greatest distance from the origin to which grains are sent during the stabilisation, while the avalanche number is the total number of topplings made. Our simulations suggest that both of these numbers have fairly simple asymptotic behaviours as functions of $\gamma$, $N$ and $M$ as $N$ tends to infinity, and we give heuristic arguments explaining how these asymptotics appear. We also provide a more detailed analysis in the case where $\gamma$ is the binomial distribution with parameter $p$, in particular when $p$ tends to $1$. We exhibit a phase transition in that limit at the scale $p \sim 1/N$.

\section{Introduction}\label{sec:intro}

\subsection{The Abelian sandpile model}\label{subsec:ASM}

The Abelian sandpile model (ASM) was originally introduced by Bak, Tang and Wiesenfeld  as an example of a system exhibiting the phenomenon known as \emph{self-organised criticality}~\cite{BTW1988}. This means that a system tunes itself to a critical state without any need for external modification of parameters~\cite{Frigg2003}. In the original ASM, grains of sand are positioned at vertices (sites) of a finite square grid. At each step, a grain of sand is added to a randomly chosen vertex on the grid. If the number of grains at a given vertex is greater than or equal to $4$, we say that the vertex has become \emph{unstable}. Unstable vertices \emph{topple}, sending one grain to each adjacent vertex on the grid. This can cause other vertices to become unstable, and these topple in turn. When a vertex at the edge of the grid topples, grains fall off the grid, exiting the system, so that this process eventually stabilises~\cite{DLT1999}. Figure~\ref{fig:ASM_example} gives an example of the evolution of the ASM. Here the colour purple is used to represent unstable vertices, while red represents grains of sand that are relocated when such a vertex topples.

\begin{figure}[h]
\centering
\includegraphics[width=0.9\textwidth]{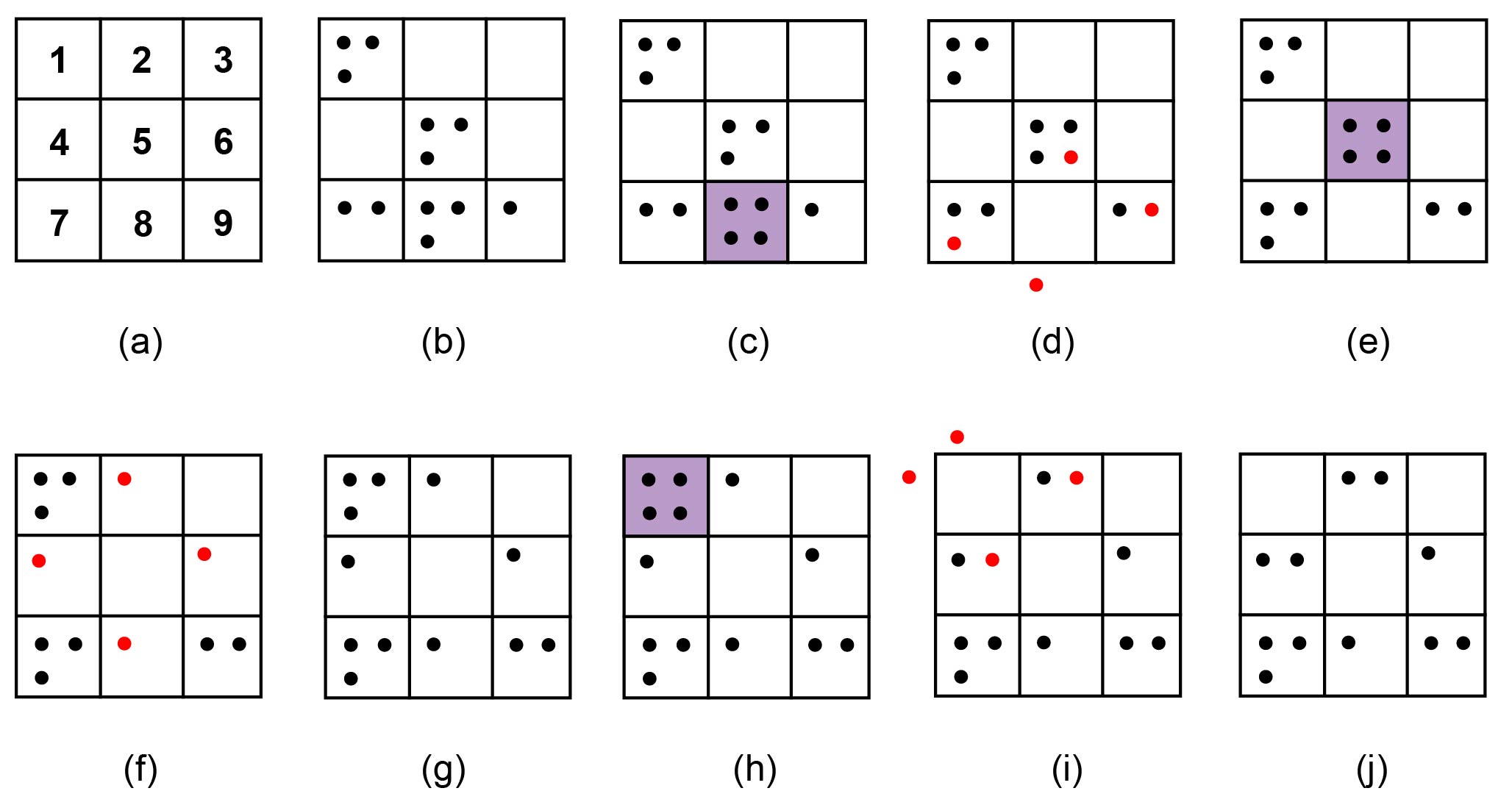}
\caption{An example of the ASM in a 3x3 grid. (a) The labels of the nine vertices, displayed for convenience; (b) The initial configuration; (c) Add one grain to vertex 8 and vertex 8 becomes unstable with four grains; (d) Vertex 8 topples once and sends one grain to each of its four neighbours;  (e) The toppling of vertex 8 causes vertex 5 to become unstable with four grains, and the grain at the bottom of the grid exits the system; (f) Vertex 5 topples once and sends one grain to each of its four neighbours; (g) The stabilised configuration after toppling once at vertex 5; (h) Add one grain to vertex 1 and vertex 1 becomes unstable with four grains; (i) Vertex 1 topples once and sends one grain to each of its four neighbours; (j) The stabilised configuration after the two grains at the top-left of the grid exit the system}
\label{fig:ASM_example}
\end{figure}

This model was formalised and extended to general graphs by Dhar~\cite{Dhar1990}, who also gave it the name of ASM. The ASM has been widely studied in Statistical Physics, Mathematics and Computer Science (under the name of ``chip-firing game'' in the latter case). It also has a number of  real-world applications such as forest-fires~\cite{RAM1999}, earthquakes~\cite{BT1989} and sediment deposits~\cite{RGF1994}. We refer the interested reader to the excellent recent survey on the topic by Klivans~\cite{Kliv2019}.

Among the various topics studied in sandpile research is the so-called \emph{single-source sandpile}. This takes place on the grid $\Z^d$. Initially, a large number $N$ of sand grains are placed at a single vertex (also called the origin). If $N \geq 2d$, the origin is unstable, and topples, yielding the toppling avalanches described above, until the process eventually stabilises. The goal is then to study this stabilisation process, usually as $N$ tends to infinity, and in particular obtain properties of the final stable configuration reached. In various simulations of this configuration, kaleidoscopic patterns will appear~\cite{B2021}, and these patterns are stated with rigourous explanations in terms of scaling limits and solutions to certain differential equations~\cite{PS2013}.

The single-source ASM has been quite extensively studied. Fey, Levine and Peres described its growth rate in~\cite{FLP2010}. In the same work, the authors gave a conjecture of dimensional reduction, which was proved by Bou-Rabee on the hypercube~\cite{B2022}. In~\cite{BG2007}, a technical tool was introduced to deal with the proof of discrete quasiconcavity for the odometer function in the 2-dimensional case ($d=2$). In~\cite{SD2012}, Sadhu and Dhar exhibited the formation of patterns for a slight variant of the single-source ASM, in which $N$ grains are dropped at a single vertex into an initial configuration satisfying some periodicity conditions.

\subsection{Stochastic variants of the ASM}\label{subsec:SSM}

In the ASM, the only level of randomness lies in the locations at which grains of sand are added to the system, after which the stabilisation process is deterministic. In particular, the single-source ASM is entirely deterministic. This may make it hard to match the complexity and apparent randomness of real-world situations. Therefore, to increase the randomness of the model, stochastic variations of the ASM have been introduced, in which the topplings are also made random. There are a number of such stochastic sandpile models in the literature.
\begin{itemize}
\item In~\cite{Manna1991}, there are two different types of grain, which cannot occupy the same vertex. When they do, one of the grains is moved to a randomly chosen neighbouring vertex instead.
\item In~\cite{DS2009}, unstable vertices lose all their grains, which are re-distributed at random to their neighbours: some neighbours may receive more than one grain, others none, while some grains may exit the system directly.
\item In~\cite{CMS2013}, unstable vertices flip a coin for each neighbour to decide which neighbours to send grains to. That is, all neighbours independently of each other receive a grain with probability $p  \in (0,1)$ (with probability $(1-p)$ that grain is kept by the toppling vertex).
\item The model in~\cite{Nunzi2016} essentially generalises the two previous models in~\cite{CMS2013, DS2009}.
\item In~\cite{KW2020}, Kim and Wang set the toppling threshold of a vertex to a (fixed) multiple $M$ of its degree. For each toppling, a random number $k \in \{1,\cdots,M\}$ is chosen according to some probability distribution $\gamma$, and each neighbour of the toppling vertex receives the same (random) number $k$ of grains. This particular SSM has been linked to the LLL algorithm for lattice-basis reduction~\cite{DKTWY2021}.
\end{itemize}

In general, the research for these stochastic variants is quite sparse compared to that of the ASM. This paper aims to remedy some of this by analysing the 2-dimensional single-source sandpile for the Kim-Wang model described above~\cite{KW2020}. To our knowledge, this is the first study of a stochastic single-source sandpile model.

Our paper is organised as follows. In the rest of this section we give a formal description of the Kim-Wang stochastic sandpile model and of the radius and avalanche numbers studied in this paper (Section~\ref{subsec:KW}), and introduce the random toppling distributions $\gamma$ considered (Section~\ref{subsec:random_multiplicity}). Section~\ref{sec:general_behav} describes the general asymptotic behaviour of the radius and avalanche numbers for these various distributions $\gamma$. We observe that the radius number behaves like $c_R(\gamma) \cdot \sqrt{\frac{N}{M}}$, while the avalanche number behaves like $c_{Av}(\gamma) \cdot \frac{N^2}{M \cdot E(\gamma)}$, where $c_R$ and $c_{Av}$ are constants that only depend on the distribution $\gamma$ (thus independent of $N$ and $M$). We also give heuristic arguments for these behaviours. Section~\ref{sec:bin_analysis} focuses on the case where $\gamma$ is a Binomial distribution with parameters $M$ and $p$ (for some $p \in (0, 1)$). In particular, we analyse in detail the limit $p \rightarrow 1$, and exhibit a phase transition when $p \sim 1/N$. Section~\ref{sec:conclusion_future} summarises our main results, and gives some possible directions for future work.

\subsection{The Kim-Wang model}\label{subsec:KW}

In this paper, we study a slight variant of the Kim-Wang model (KW) on the grid $\Z^2$. The reason for our choice of model is that the topplings of KW are symmetric (when a vertex topples it sends the same number of grains to each of its neighbours), thus emulating some of the symmetries of $\Z^2$. We now describe the model in more detail, including formal definitions. Fix some positive integer $M$, called the \emph{multiplicity} of the model. Also fix some probability distribution $\gamma$ whose support is a subset of $\{1, \cdots, M\}$ (in practice it will usually be equal). A \emph{configuration} is a vector $(c(v))_{v \in \Z^2}$ such that for all $v \in \Z^2$, $c(v)$ is a non-negative integer, and $c$ has finite support (i.e. $c(v) = 0$ for all $v$ outside a finite subset of $\Z^2$). For $v \in \Z^2$, we denote by $\ind_v$ the configuration such that $\ind_v(v) = 1$ and $\ind_v(w) = 0$ for $w \neq v$.

For a configuration $c$, we say that a vertex $v \in \Z^2$ is \emph{stable} if $c(v) < 4M$, and unstable otherwise. Unstable vertices \emph{topple}, by sending $k$ grains of sand to each of their four neighbours in $\Z^2$ (the four neighbours of $v = (x,y) \in \Z^2$ are $(x \pm 1, y)$ and $(x, y \pm 1)$), where $k$ is a random number drawn according to the distribution $\gamma$. Formally, for $v \in \Z^2$ and $i > 0$, we define $\Topp{\gamma}{v, i}$, the $i$-th \emph{toppling operator} at $v$, by:
\begin{equation}\label{eq:toppling_op}
\Topp{\gamma}{v, i}(c) = c - 4X_{v,i} \ind_v + \sum\limits_{w \sim v} X_{v,i} \ind_w,
\end{equation}
where the sum is over all neighbours $w$ of $v$, and the $\left(X_{v,i}\right)_{v \in \Z^2,\, i \in \N}$ are i.i.d. random variables with distribution $\gamma$. Note that the total number of grains $\sum\limits_{v \in \Z^2} c(v)$ is preserved through this operation.

The single-source SSM is then defined as follows. Start from an initial configuration $c = c^0 := N \cdot \ind_{(0,0)}$ for some positive integer $N$. If $v$ in $c$ is unstable and has previously toppled $i-1$ times, set $c = \Topp{\gamma}{v, i}(c)$. 
Repeat the process until no more vertices can be toppled, i.e. $c$ has become stable, and set $c'$ to be the stable configuration reached. For $v \in \Z^2$, define the \emph{toppling number} $N(v)$ as the number of times $v$ topples in the stabilisation $c^0 \rightsquigarrow c'$. It is possible to show (see e.g.~\cite{Nunzi2016}) that the toppling numbers $N(v)$ do not depend on the order in which unstable vertices are toppled, and neither does the final configuration $c'$ reached, which we denote by $\Stab^{(\gamma; N, M)}$. 

We now define the two main objects of our study. The \emph{radius number} $R(\gamma; N, M)$ is the max distance from the origin where grains are sent in the stabilisation above, that is:
\begin{equation}\label{eq:def_radius}
R(\gamma; N, M) = \max\limits_{v \in \Z^2} \{ \| v \|_2; \, \Stab^{(\gamma; N, M)}(v) > 0 \},
\end{equation}
where $\| \cdot \|_2$ denotes the usual L2 norm in $\Z^2$. The \emph{avalanche number} $Av(\gamma; N, M)$ is the total number of topplings in the stabilisation $c^0 \rightsquigarrow \Stab^{(\gamma; N, M)}$, i.e.:
\begin{equation}\label{eq:def_avalanche}
Av(\gamma; N, M) = \sum\limits_{v \in \Z^2} N(v),
\end{equation}
where the $N(v)$ are the toppling numbers defined above. Since the configuration $\Stab^{(\gamma; N, M)}$ and the toppling numbers $N(v)$ do not depend on the toppling order, neither do the radius and avalanche numbers.

\subsection{The random multiplicity distributions}\label{subsec:random_multiplicity}

We refer to the distribution $\gamma$ introduced in Section~\ref{subsec:KW} as the \emph{random multiplicity} of the model, since it is used to generate the multiplicity of each toppling. 
Recall that each time a vertex becomes unstable, we re-sample $\gamma$, so that the number of grains toppled could be different each time. 
Because the only condition on $\gamma$ is that its support be included in $\{1, \cdots , M\}$, the selection of $\gamma$ could be flexible, with a number of possible distributions. In this paper, we study the single-source SSM for six different $\gamma$ distributions, given as follows.
\begin{itemize}
\item \textbf{Deterministic distribution:} $\gamma = M$ a.s., so that each time a vertex topples, it sends $M$ grains to each of its neighbours. This is equivalent to the ASM, up to rescaling (see Section~\ref{subsec:heuristics_asymptotics}).

\item \textbf{Uniform distribution:} $\gamma$ is a uniform random variable with range from $1$ to $M$, i.e. $P(\gamma = k) = \frac{1}{M}$ for all $k \in \{1, \cdots , M\}$. 

\item \textbf{Binomial distribution:} $\gamma$ is a binomial random variable with parameters $M, p$, for some $p \in (0, 1]$, conditioned to be non-zero (so that $\gamma \geq 1$ a.s.). That is, $P(\gamma = k) = \dfrac{\binom{M}{k} p^k (1-p)^{M-k}}{1 - (1-p)^M}$ for all  $k \in \{1, \cdots , M\}$. 

\item \textbf{Log-Law distribution:} this is the probability distribution such that its cumulative distribution function at $k$ is proportional to $\log(k+1)$ for $k \geq 1$. In order to get the desired support for $\gamma$, the correct renormalisation yields $P(\gamma = k) = \log_{(M+1)}{(k+1)}-\log_{(M+1)}{k}$ for all $k \in \{1, \cdots , M\}$, where $\log_a(b) = \frac{\log(b)}{\log(a)}$. 
This gives the desired cumulative distribution function $F(k):= P(\gamma \leq k) =\log_{(M+1)}{(k+1)}$.

\item \textbf{Power-Law distribution:} $\gamma$ follows a finite-support power-law distribution with parameter $s > 0$. That is, $P(\gamma=k)=\frac{k^{-s}}{Z(s,M)}$ for all $k \in \{1, \cdots , M\}$, where $Z(s, M) := \sum_{j=1}^{M} j^{-s}$ is the corresponding partition function.

\item \textbf{Always-1 distribution:} $\gamma = 1$ a.s., so that each time a vertex topples, it sends one grain to each of its neighbours. Note that when $M=1$, the Always-1 distribution and Deterministic distribution are the same, with the same toppling rules as the standard ASM. In general, however, the Always-1 distribution differs from the ASM due to its toppling threshold, which is still $4M$ (as opposed to $4$ in the ASM). This distribution is introduced in order to study the limit of the Binomial distribution as $p$ tends to 0. (see Section~\ref{subsec:p_to_0} for more details).

\end{itemize}

\section{General behaviour}\label{sec:general_behav}

Recall briefly the goal of our study. We start with $N$ grains at the origin of $\Z^2$, and study the stabilisation process as $N$ tends to infinity. The stabilisation here is the random process of the KW model: the toppling threshold at each vertex is $4M$ for some positive integer $M$, and topplings send a random number $k$ of grains to each neighbour according to some probability distribution $\gamma$ on $\{1, \cdots, M\}$. Our research focuses on two global parameters of the system:
\begin{itemize}
\item the radius number $R(\gamma; N, M)$, describing how far away from the origin this process reaches;
\item the avalanche number $Av(\gamma; N, M)$, describing the total number of topplings that occur in this process.
\end{itemize}
These parameters are random variables, whose distribution depends on the distribution $\gamma$, as well as the parameters $N$ and $M$. We use the term “global” here in contrast to ``local'' parameters such as the number of times each individual vertex topples (called the \emph{odometer function} in sandpile research, see e.g. \cite{PS2013}). Such local phenomena may also be of interest in future study (see Section~\ref{sec:conclusion_future} for more details on this).

\subsection{The radius and avalanche constants}\label{subsec:radius_avalanche_const}

For the model developed in this paper, we observe the following general asymptotic behaviour of the radius number, as $N$ tends to infinity with $M$ fixed:
\begin{equation}
R(\gamma; N, M) \sim c_R(\gamma) \cdot \sqrt{\frac{N}{M}},\label{eq:asympt_radius}
\end{equation}
where $c_{R}(\gamma)$ is a constant whose value only depends on $\gamma$ and neither on $N$ nor on $M$.

For the avalanche number, the observed asymptotic behaviour is:
\begin{equation}
Av(\gamma; N, M) \sim c_{Av}(\gamma) \cdot \frac{N^2}{(M \cdot E(\gamma))},\label{eq:asympt_avalanche}
\end{equation}
where $E(\gamma)$ is the expected value of the $\gamma$ distribution, and again $c_{Av}(\gamma)$ is a constant whose value only depends on $\gamma$ and neither on $N$ nor on $M$. We will refer to $c_R$, resp. $c_{Av}$, as the \emph{radius constant}, resp. \emph{avalanche constant}, of the model.

Here Equations~\eqref{eq:asympt_radius} and \eqref{eq:asympt_avalanche} are what we observe on the data, and in Section~\ref{subsec:heuristics_asymptotics} we will give heuristics, i.e. non-rigourous arguments, for why these behaviours seem plausible. Section~\ref{subsec:numeric_estimates} is then dedicated to showing the data, and giving numerical estimates for the radius and avalanche constants.

We end this section with a short discussion on the definition of the radius number in Equation~\eqref{eq:def_radius}. In this, we have chosen L2 as our representative distance. It may seem \textit{a priori} that L1 would have been a more appropriate choice, since it appears to better account for the propagation of avalanches (the neighbours of a vertex are those at L1-distance 1). 

However, Figure~\ref{fig:shapes} below shows the shape of the final stable configuration reached $\Stab^{(\gamma; N, M)}$ for different distributions $\gamma$ and values of $M, N$. Here, we use a sliding RGB colour scale for the number of grains at each vertex. The colour blue is used for vertices with $1$ grain, green for vertices with $2M$ grains, and red for vertices with $4M - 1$ grains. Linear interpolation is then applied on the RGB scale from blue to green (number of grains from $1$ to $2M$ and green to red (number of grains from $2M$ to $4M-1$). We can see from these simulations that the shapes are circular, or near-circular. This suggests that the L2 distance (or Euclidean distance) we use is indeed the more appropriate choice in the calculation of the radius number. Some brief further observations on these shapes are made in Section~\ref{subsec:future}.

\begin{figure}[h]
\centering
\includegraphics[width=0.9\textwidth]{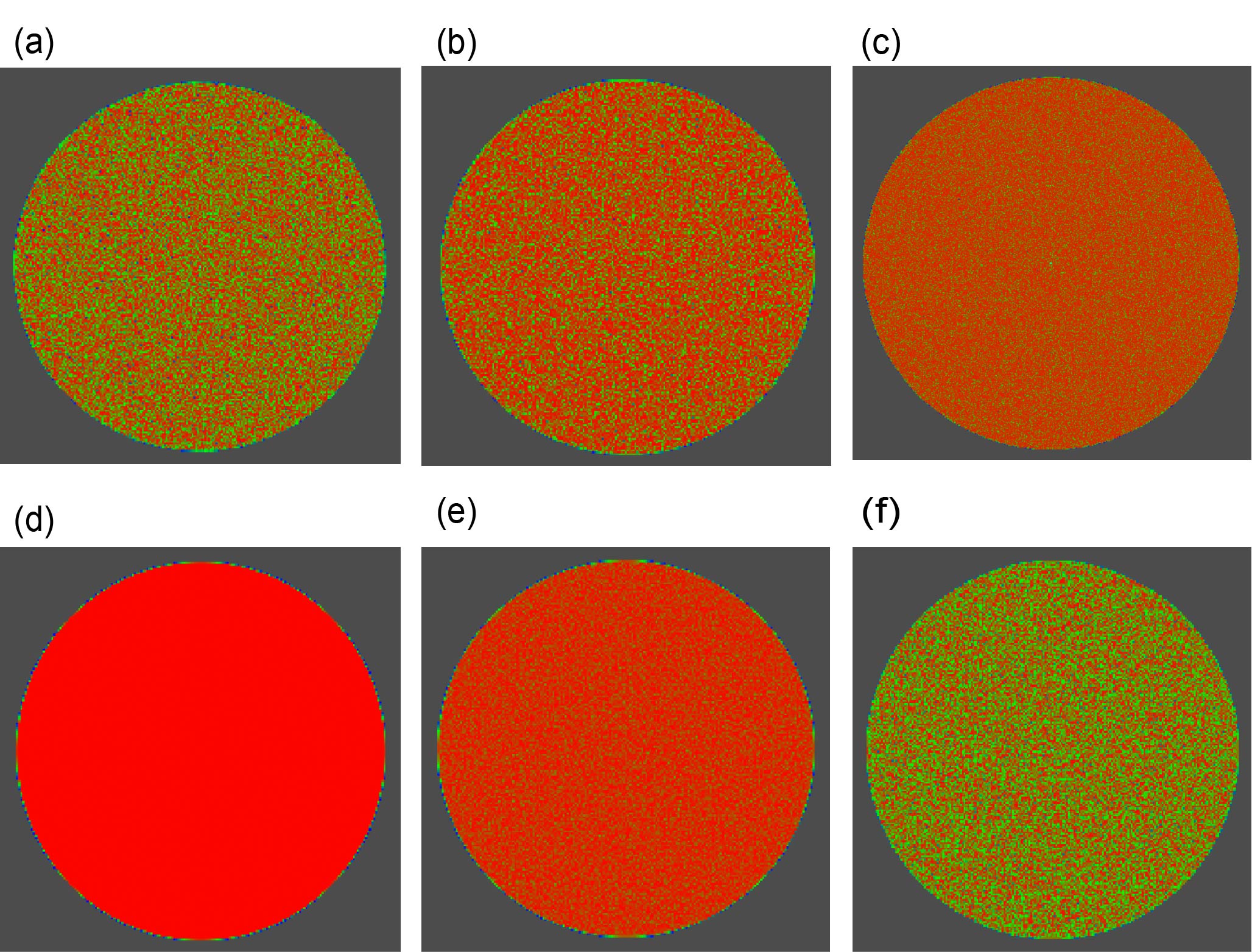}
\caption{Some shapes after stabilisation. (a) Uniform distribution for $M=100$, $N=10^7$; 
(b) Log-Law distribution for $M=100$, $N=10^7$; 
(c) Power-Law distribution for $M=10$, $N=9\cdot 10^6$, and $s=2$;
(d) Always-1 distribution for $M=100$, $N=10^7$;
(e) Binomial distribution for $M=100$, $N=10^7$ and $p=0.25$;
(f) Binomial distribution for $M=100$, $N=10^7$ and $p=0.75$}
\label{fig:shapes}
\end{figure}

\subsection{Heuristics for the asymptotics}\label{subsec:heuristics_asymptotics}

In this section, we are going to give informal arguments for why the radius and avalanche numbers behave as described in Equations~\eqref{eq:asympt_radius} and \eqref{eq:asympt_avalanche} from the previous section. First, we claim that if $N$ is a multiple of $M$, we have $Av(\Det; N, M) = Av(\Det; N/M, 1)$.

To see this, take two small examples with $N=40$ and $M=2$ for $\Det(N, M)$ and $\Det(N/M, 1)$, which are shown in Figures~\ref{fig:simulation_Det(N, M)} and \ref{fig:simulation_Det(N/M, 1)} respectively. Here the notation $\Det(N, M)$ means the single-source Deterministic model (from Section~\ref{subsec:random_multiplicity}) with $N=40$ grains at the initial central vertex and multiplicity $M=2$, meaning that topplings occur when a vertex has at least $4M = 8$ grains, sending $M = 2$ grains to each neighbour. $\Det(N/M, 1)$ means the single-source Deterministic model with $N/M=20$ grains at the initial central vertex and multiplicity $1$, i.e. topplings occur when a vertex has at least $4$ grains, sending one grain to each neighbour. Note that $\Det(k, 1)$ is in fact the single-source ASM with $k$ initial grains at the central vertex.

In Figure~\ref{fig:simulation_Det(N, M)}, the processes (a)-(d) describe the stabilisation of  $\Det(N, M)$, and in Figure~\ref{fig:simulation_Det(N/M, 1)}, the processes (a)-(d) show the stabilisation of $\Det(N/M, 1)$. For convenience, the central vertex is denoted in bold and is surrounded by black bold lines, while unstable vertices are represented in italic.

\begin{figure}[h]
\centering
\includegraphics[width=0.9\textwidth]{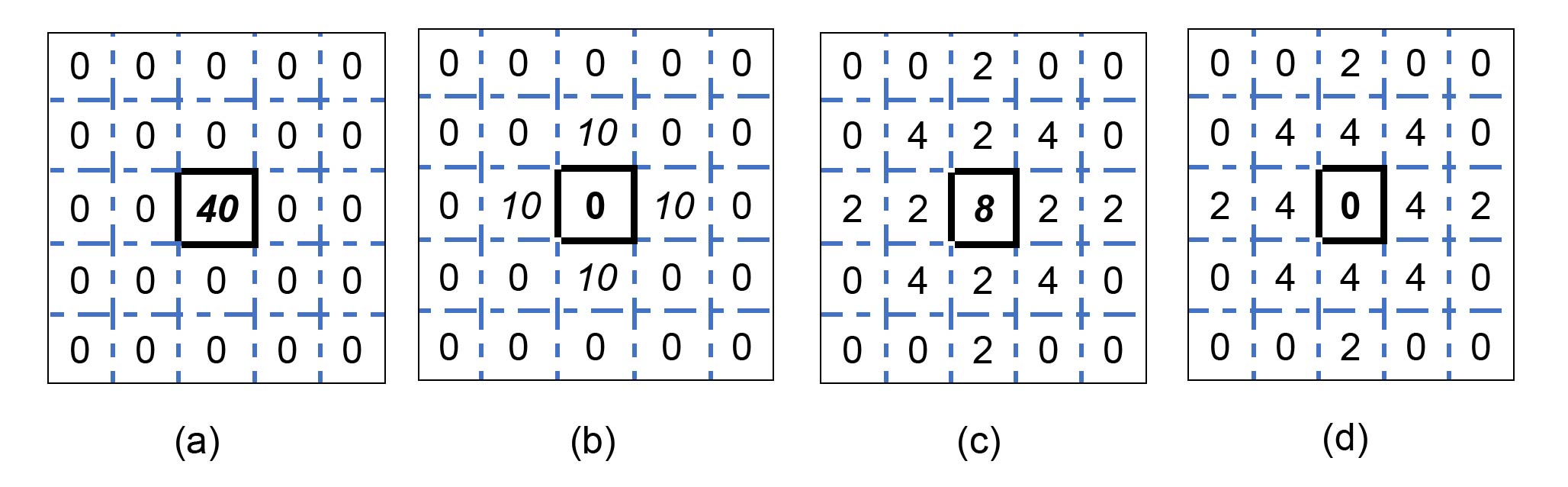}
\caption{Simulation for the stabilisation of $\Det(N, M)$ in a 5x5 grid. (a) The initial configuration for $N=40$, $M=2$; (b) The configuration after five topplings of the central vertex; (c) The configuration after one toppling for the four unstable neighbours of the central vertex; (d) The stabilised configuration after one further toppling of the central vertex}
\label{fig:simulation_Det(N, M)}
\end{figure}

For $\Det(N, M)$ shown in Figure~\ref{fig:simulation_Det(N, M)}, we start with 40 grains at the central vertex in (a). Recall that the toppling threshold here is $4M = 8$, so this vertex is unstable. In fact, we can topple this vertex five times successively (because the number of topplings needs to reduce its number of grains to be less than the toppling threshold 8), after which the four neighbours have 10 grains each (they pick up $M=2$ grains in each toppling), giving us the result in (b). After (b), these four neighbours are in turn unstable, and they will each topple once, sending two grains to each neighbour, resulting in (c). Finally, in (c), the central vertex is unstable and will topple again, and the final stable configuration shown in (d) will be reached. In this stabilisation, the radius number $R(\Det; N, M)$ is $2$, and the avalanche number $Av(\Det; N, M)$ is $5 + 4 + 1 = 10$.

\begin{figure}[h]
\centering
\includegraphics[width=0.9\textwidth]{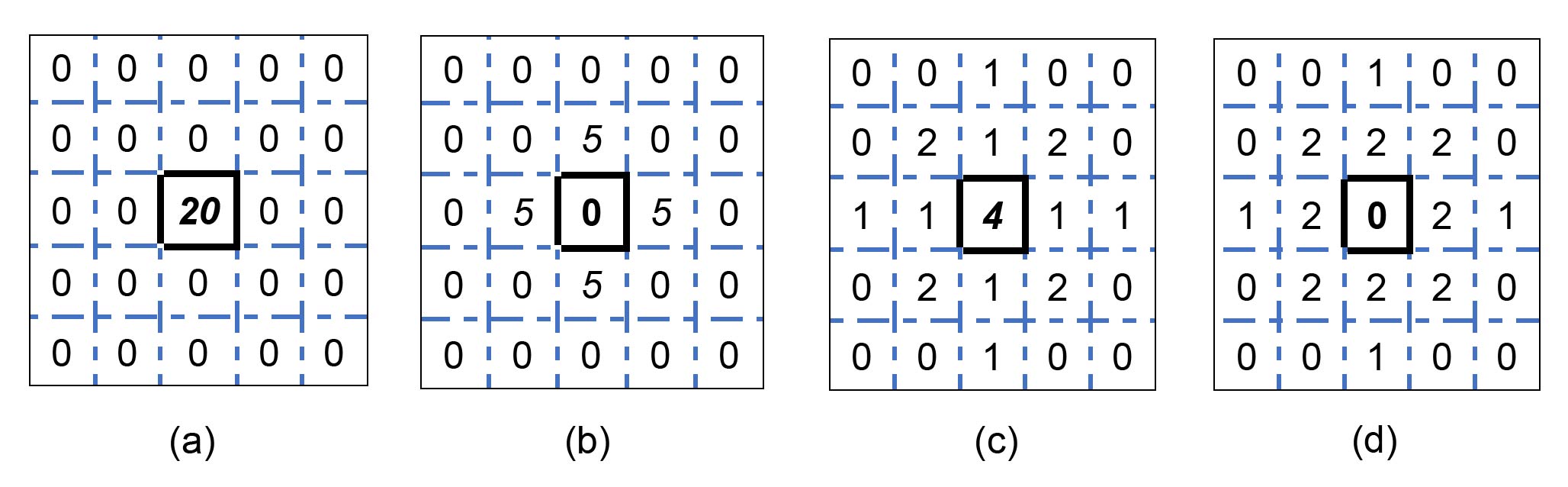}
\caption{Simulation for the stabilisation of $\Det(N/M, 1)$ in a 5x5 grid. (a) The initial configuration for $N=20$, $M=1$; (b) The configuration after five topplings of the central vertex; (c) The configuration after one toppling for the four unstable neighbours of the central vertex; (d) The stabilised configuration after one further toppling of the central vertex}
\label{fig:simulation_Det(N/M, 1)}
\end{figure}

Now do the same process again to check the behaviour of $\Det(N/M, 1)$ by observing the demonstrated processes (a)-(d) of Figure~\ref{fig:simulation_Det(N/M, 1)}. This time, we start in (a) with $20$ grains initially, and multiplicity $1$ (so the toppling threshold is $4 \cdot 1 = 4$). After toppling the unstable central vertex in (a) five times, the situation in (b) is displayed (again, we keep toppling the central vertex until its number of remaining grains is less than the toppling threshold, which is $4$ here). In (b), the four neighbours of the central vertex have become unstable. Topple these each once to obtain the result in (c). Finally, in (c), topple the central vertex again, to get the stable configuration reached in (d). In this stabilisation, the radius number $R(\Det; N/M, 1)$ is 2, and the avalanche number $Av(\Det; N/M, 1)$ is $5 + 4 + 1 = 10$.

Note that the toppling sequence is exactly the same for the $\Det(N, M)$ and $\Det(N/M, 1)$ cases, which implies immediately that the radius and avalanche numbers will also be the same. This also means that the final configuration that we have reached for $\Det(N, M)$ is just $M$ times the final configuration reached for $\Det(N/M, 1)$. This last observation is not quite true if $N$ is not a multiple of $M$ (e.g. we may get some sort of “remainder” at the central origin), but any error terms are relatively small (of size at most the toppling threshold $4M$).

From the above, we have that $R(\Det; N, M) = R(\Det; N/M, 1)$. In addition, because $\Det(k, 1)$ is just the standard ASM, from results in~\cite{SD2010} and references therein, we know that $R(\Det; k, 1) \sim c \cdot \sqrt{k}$ as $k$ tends to infinity (this is the case of no sink sites). Hence, for fixed $M$, we get:
\begin{equation}\label{eq:radius_number_behave}
R(\Det;N,M) \approx c \cdot \sqrt{\frac{N}{M}},
\end{equation}
as $N$ tends to infinity. In the case of the radius number, the data suggests that the behaviour for general gamma distributions is the same as that of the Deterministic model above, with only the multiplicative constant $c$ changing. This gives the claimed behaviour and definition of the radius constant from Equation~\eqref{eq:asympt_radius}.

As well as the behaviour of the radius number, from above, we have obtained $Av(\Det; N, M) = Av(\Det; N/M, 1)$. As a consequence of results in \cite{PS2013}, we know that $Av(\Det; k, 1)$ behaves like $c \cdot k^2$ as $k$ tends to infinity. Therefore, for fixed $M$, we get:
\begin{equation} \label{eq:avalanche_number_behav}
Av(\Det;N,M) \approx c \cdot \left(\frac{N}{M}\right)^2,
\end{equation}
as $N$ tends to infinity. In this case, the data suggests the behaviour is slightly more complex in the SSM. We now give some heuristic (informal) arguments for the claimed behaviour from Equation~\eqref{eq:asympt_avalanche}.

Consider a toppling sequence of vertices denoted by $S = v_1, \cdots , v_k$. Define the toppling sequence $S^M := v_1, \cdots , v_1, \cdots , v_k, \cdots , v_k$, where each vertex of the sequence $S$ appears $M$ times successively. It is clear that applying the toppling sequence $S$ in the $\Det(N, M)$ model will yield the exact same configuration as applying the toppling sequence $S^M$ in $\A1(N, M)$, since toppling a vertex once in $\Det(N, M)$ has the same effect as toppling a vertex $M$ times in $\A1(N,M)$. Here, $\A1(N,M)$ denotes the single-source sandpile for  the Always-1 distribution with $N$ grains and multiplicity $M$.

This means that if $S$ is the stabilisation sequence of $\Det(N, M)$, then applying $S^M$ in $\A1(N, M)$ will also give a stable configuration (the toppling threshold is the same in both models). By minimality of the stabilisation sequence (see e.g.,~\cite{PS2013}), this implies
\begin{equation}
Av(\A1;N,M)\leq M \cdot Av(\Det;N,M). \label{eq:avalanche_Always-1_Det}
\end{equation}
In fact, this inequality is true locally at each vertex. That is, for any vertex $v \in \Z^2$, the number of times $v$ topples in the stabilisation of $\A1(N, M)$ is bounded above by $M$ times the number of times $v$ topples in the stabilisation of $\Det(N, M)$.

Now fix some integer $k$ such that $M/k$ is also an integer. Toppling a vertex once in $\Det(N, M)$ has the same effect as toppling a vertex $M/k$ times in $\A k(N, M)$, where $\A k(N, M)$ is the single-source sandpile with $N$ grains and multiplicity $M$, in which all topplings send exactly $k$ grains to each neighbour of a toppling vertex. As in the derivation of Equation~\eqref{eq:avalanche_Always-1_Det} above, this gives
\begin{equation} \label{eq:avalanche_Always-k_Det}
Av(\A k;N,M) \leq \left(\frac{M}{k}\right) \cdot Av(\Det;N,M).
\end{equation}

Now consider the case of more general probability distributions $\gamma$. We apply some form of the law of large numbers to the random topplings. Because the avalanche numbers are large as $N$ tends to infinity, it seems plausible to be able to plug $k = E(\gamma)$ into Equation~\eqref{eq:avalanche_Always-k_Det}, perhaps at the cost of a multiplicative constant $c_1 > 1$ in the right-hand side, yielding:
\begin{equation} \label{eq:avalanche_upperbound}
Av(\gamma;N,M) \leq c_1 \cdot \left(\frac{M}{E(\gamma)}\right) \cdot Av(\Det;N,M). 
\end{equation}

Equation~\eqref{eq:avalanche_upperbound} only gives an upper-bound for the scale. Nevertheless, it seems likely that in Equation~\eqref{eq:avalanche_Always-k_Det} the bound is relatively tight. Recall here that the difference is given by the gap in size between the stabilisation sequence for $\A k(N, M)$ and the stabilisation sequence of $\Det(N, M)$ applied $\left( \frac{M}{k} \right)$ times. In other words, it measures how close the latter sequence is to being the minimal stabilising sequence. It seems plausible that this gap is small, yielding at most a linear error term. This would give the following inequality:
\begin{equation} \label{eq:avalanche_Always-k_greaterthan_Det}
Av(\A k;N,M) \geq c \cdot \left(\frac{M}{k}\right) \cdot Av(\Det;N,M),
\end{equation}
for some constant $c < 1$.

From Equation~\eqref{eq:avalanche_Always-k_greaterthan_Det}, a similar use of the law of large numbers would then yield:
\begin{equation} \label{eq:avalanche_lowerbound}
Av(\gamma;N,M) \geq c_2 \cdot \left(\frac{M}{E(\gamma)}\right) \cdot Av(\Det;N,M),
\end{equation}
for some constant $c_2 < 1$.

Combining Equations~\eqref{eq:avalanche_upperbound} and \eqref{eq:avalanche_lowerbound} for $Av(\gamma; N, M)$ with Equation~\eqref{eq:avalanche_number_behav} which describes the behaviour of $Av(\Det; N, M)$ gives
\begin{equation} \label{eq:avalanche_range}
c_2 \cdot \frac{N^2}{M \cdot E(\gamma)} \leq Av(\gamma;N,M) \leq c_1 \cdot \frac{N^2}{M \cdot E(\gamma)},
\end{equation}
which shows that $Av(\gamma; N, M)$ should indeed behave like $c_{Av}(\gamma)\cdot\frac{N^2}{M \cdot E(\gamma)}$, as in Equation~\eqref{eq:def_avalanche}.

\subsection{Numeric estimates for avalanche and radius constants}\label{subsec:numeric_estimates}

In this section we present our experimental data that gives numerical estimates for the radius and avalanche constants as defined in Equations~\eqref{eq:def_radius} and \eqref{eq:def_avalanche}. To estimate the radius constant, we plot the quotient $\dfrac{R(\gamma; N, M)}{\sqrt{\frac{N}{M}}}$ as a function of $N$ for the various $\gamma$ distributions from Section~\ref{subsec:random_multiplicity} and a few different values of $M$.
For the avalanche constant, we plot the quotient $\dfrac{Av(\gamma; N, M)}{\frac{N^2}{M \cdot E(\gamma)}}$.

For the expected value of $\gamma$, in some cases we have explicit formulae, while in others (Low-Law and Power-Law) no such formula exists. In the latter case we simply use numerical values derived from the formula $E(\gamma) = \sum\limits_{k = 1}^M k \cdot P(\gamma = k)$. Table~\ref{table:expected_values} summarises the various $E(\gamma)$ values used.

\begin{table}[ht]
\begin{center}
\caption{\label{table:expected_values}The expected values for the various $\gamma$ distributions}
\begin{tabular}{cc}
\toprule
$\gamma$ distribution & $E(\gamma)$ \\
\midrule
$\Det(N, M)$ & $M$ \\[5pt]
$\A1(N,M)$ & $1$ \\[5pt]
Uniform & $\dfrac{M+1}{2}$ \\[10pt]
Binomial & $\dfrac{M \cdot p}{1 - (1-p)^M}$ \\[11pt]
Low-Law & Numerical value \\[4pt]
Power-Law & Numerical value \\
\bottomrule
\end{tabular}
\end{center}
\end{table}

We plot the two quotient estimates above for the six distributions of Table~\ref{table:expected_values} as functions of $N$. Our data calculates these quotients for $M=10$, $100$ and $1000$, with $N$ ranging from $10^4$ to $10^7$. In the case of the Power-Law distribution, we limit our study to the case $s=2$, as simulations were quite slow for that particular model (see Section~\ref{sec:conclusion_future} for a discussion on this). These plots are shown in Figures~\ref{fig:radius_avalanche_Det}--\ref{fig:radius_avalanche_PowerLaw}. We do not plot the estimates in the Binomial distribution case in this section, as that will be the focus of a more detailed analysis in Section~\ref{sec:bin_analysis}.

\begin{figure}[ht]
\centering
\includegraphics[width=0.9\textwidth]{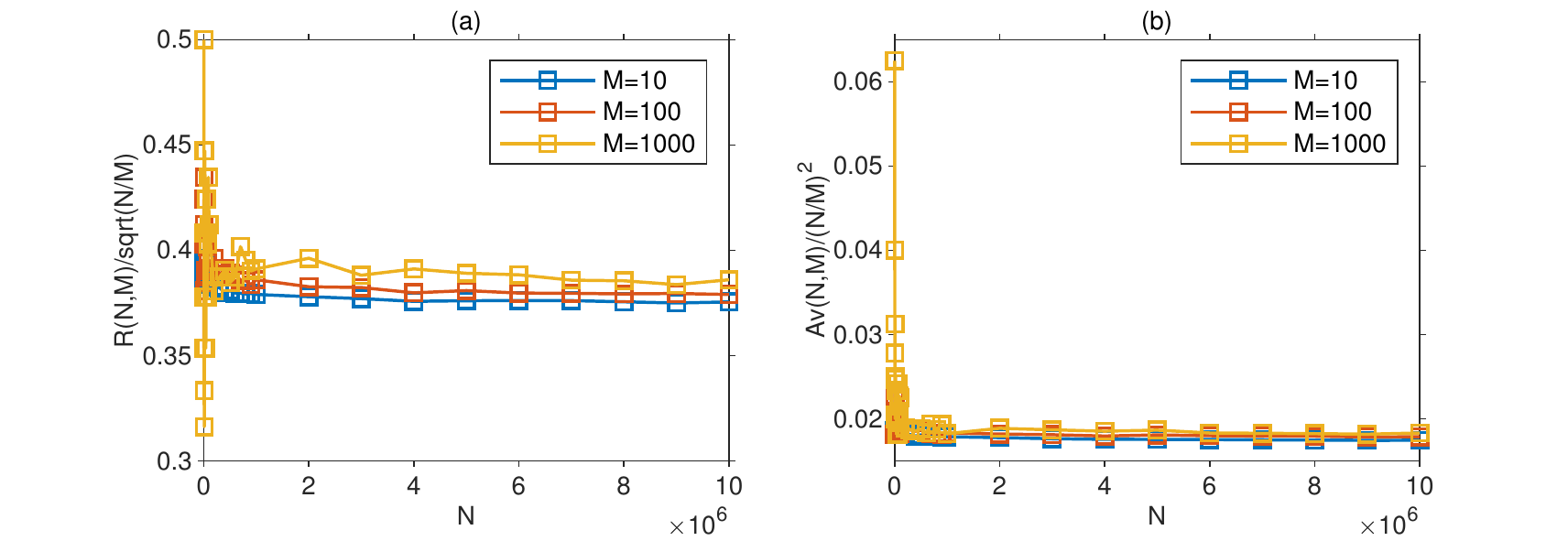}
\caption{Estimates for the radius constant (a) and avalanche constant (b) for the Deterministic distribution}
\label{fig:radius_avalanche_Det}
\end{figure}

\begin{figure}[ht]
\centering
\includegraphics[width=0.9\textwidth]{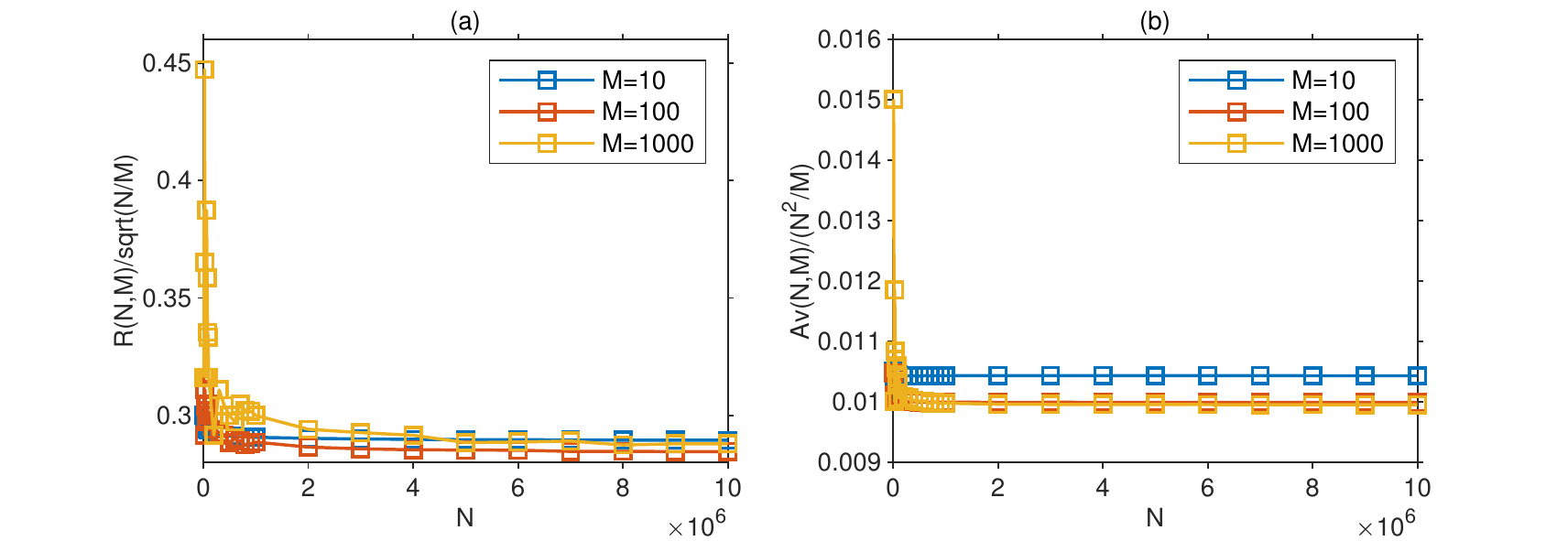}
\caption{Estimates for the radius constant (a) and avalanche constant (b) for the Always-1 distribution}
\label{fig:radius_avalanche_Always-1}
\end{figure}

\begin{figure}[ht]
\centering
\includegraphics[width=0.9\textwidth]{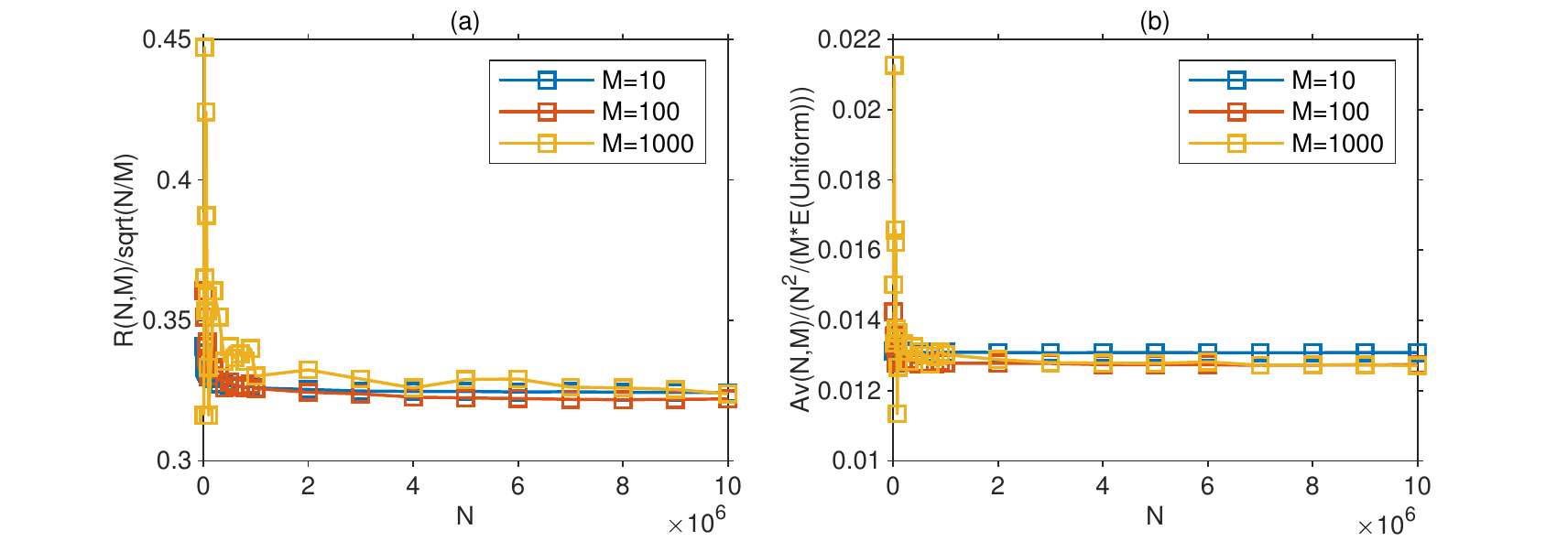}
\caption{Estimates for the radius constant (a) and avalanche constant (b) for the Uniform distribution}
\label{fig:radius_avalanche_Uniform}
\end{figure}

\begin{figure}[ht]
\centering
\includegraphics[width=0.9\textwidth]{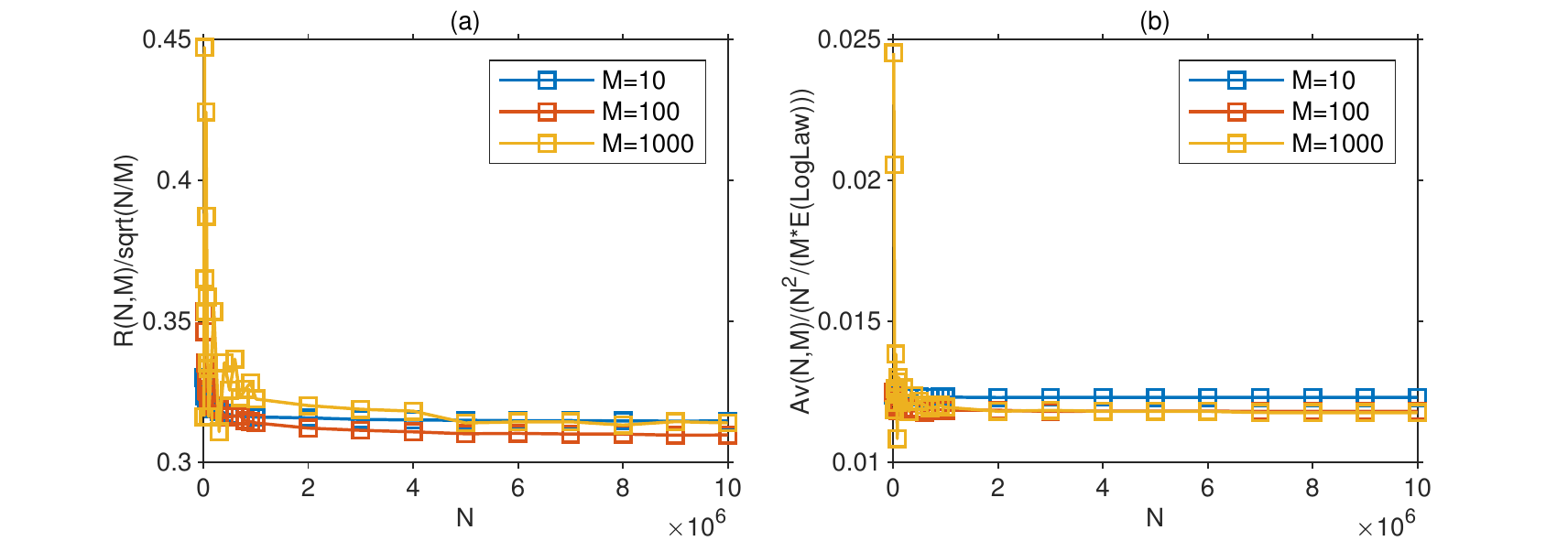}
\caption{Estimates for the radius constant (a) and avalanche constant (b) for the Log-Law distribution}
\label{fig:radius_avalanche_LogLaw}
\end{figure}

\begin{figure}[ht]
\centering
\includegraphics[width=0.9\textwidth]{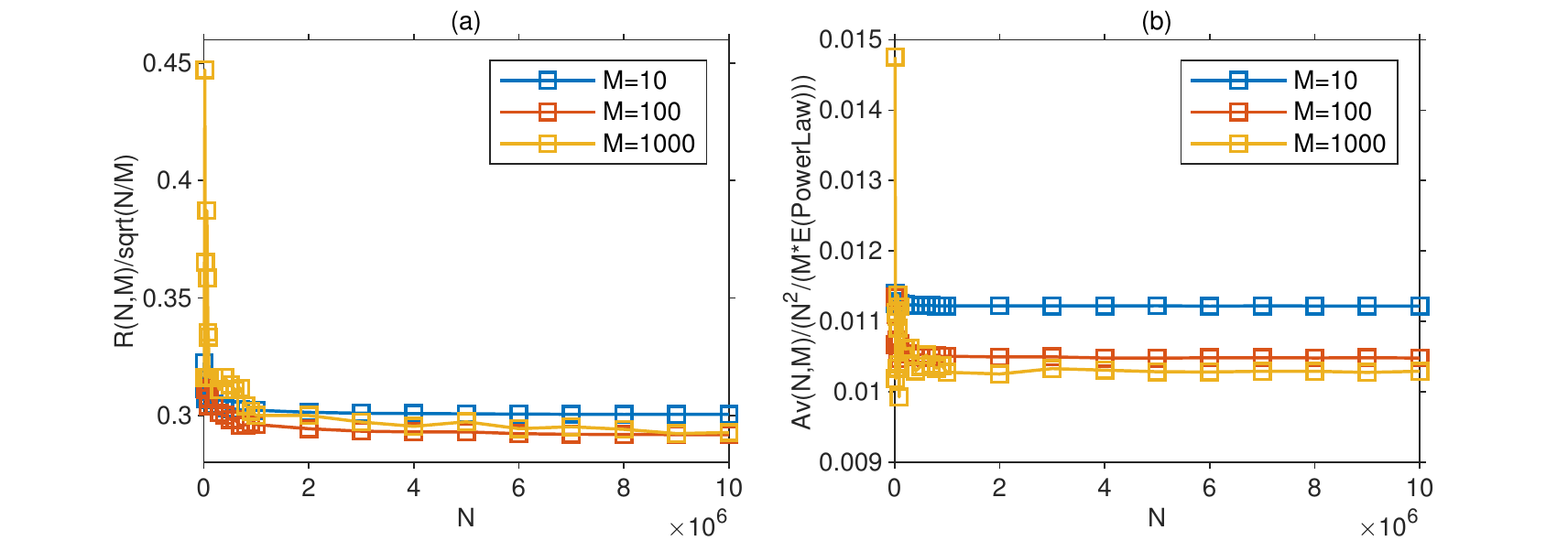}
\caption{Estimates for the radius constant (a) and avalanche constant (b) for the Power-Law distribution (here $s=2$)}
\label{fig:radius_avalanche_PowerLaw}
\end{figure}

We observe the claimed behaviour from Equations~\eqref{eq:def_radius} and \eqref{eq:def_avalanche}, namely that the quotients $\dfrac{R(\gamma; N, M)}{\sqrt{\frac{N}{M}}}$ and $\dfrac{Av(\gamma; N, M)}{\frac{N^2}{M \cdot E(\gamma)}}$ converge to some limits that do not depend on $M$ as $N$ tends to infinity. Table~\ref{table:five_radius_avalanche} shows the numerical estimates obtained for both of these quotient limits, called radius and avalance constants respectively.

\begin{table}[ht]
\begin{center}
\begin{minipage}{330pt}
\caption{\label{table:five_radius_avalanche}Numerical estimates for the radius constant and avalanche constant for the Deterministic, Always-1, Uniform, Log-Law and Power-Law distributions}
\begin{tabular}{@{}lll@{}}
\toprule 
Distribution type & The radius constant & The avalanche constant\\
\midrule
Deterministic & [0.37548, 0.386005] & [0.0174483, 0.0183092]\\[8pt]
Always-1 & [0.28457, 0.28943] & [0.00995538, 0.0104344]\\[8pt]
Uniform & [0.322025, 0.323883] & [0.0127099, 0.0130817]\\[8pt]
Log-Law & [0.309742, 0.314006] & [0.0117649, 0.0123048]\\[8pt]
Power-Law & [0.29177, 0.300483] & [0.0102943, 0.0112208]\\
\bottomrule
\end{tabular}
\end{minipage}
\end{center}
\end{table}

We observe the following relationships between the expected values of $\gamma$ and the two constants, suggesting that the constants are both increasing in the expected value $E(\gamma)$.
\begin{equation}\label{eq:const_inc_E}
\begin{split}
E(\Det)  & >  E(\mathrm{Uniform}) >   E(\mathrm{Log\text{-}Law}) > E(\mathrm{Power\text{-}Law})  >  E(\A1),   \\[6pt]
c_{R}(\Det) & >  c_{R}(\mathrm{Uniform})  > c_{R}(\mathrm{Log\text{-}Law}) > c_{R}(\mathrm{Power\text{-}Law}) > c_{R}(\A1), \\[6pt]
c_{Av}(\Det) & >  c_{Av}(\mathrm{Uniform})  > c_{Av}(\mathrm{Log\text{-}Law})> c_{Av}(\mathrm{Power\text{-}Law})  > c_{Av}(\A1).\\[6pt]
\end{split}
\end{equation}

\section{Binomial distribution analysis}\label{sec:bin_analysis}

In this section, we analyse in more depth the case where $\gamma$ is a Binomial distribution with parameters $M$ and $p$, conditioned to be non-zero. We start by looking at the general case where $p$ is fixed, before focusing on the limits $p \rightarrow 1$ and $p \rightarrow 0$.

\subsection{General p}\label{subsec:general_p}

Figure~\ref{fig:radius_avalanche_bin_0to1} plots the radius and avalanche constants as functions of $p$ over the interval $[0,1]$. Here we fix $M=100$ and take the values of the quotient estimates from Section~\ref{subsec:numeric_estimates} for $N=10^7$.  For $p=0$ we took the data for the Always-1 distribution here (see Section~\ref{subsec:p_to_0} for more details on why this is the appropriate choice). For $p=1$, we took the data for the Deterministic distribution, to which this case corresponds.

\begin{figure}[h]
\centering
\includegraphics[width=0.9\textwidth]{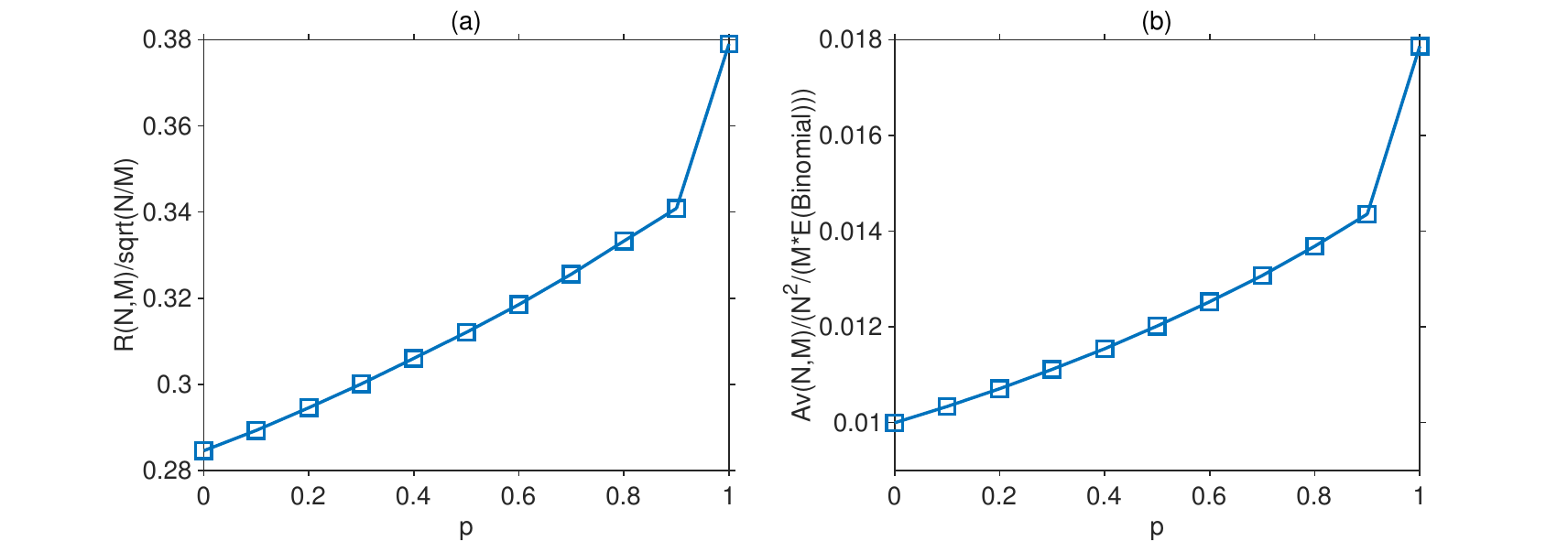}
\caption{Estimates for the radius constant (a) and avalanche constant (b) for Binomial distribution from $p=0$ to $p=1$ when $M=100$ and $N=10^7$}
\label{fig:radius_avalanche_bin_0to1}
\end{figure}

We observe two things. Firstly, in the interval $[0, 0.9]$ the two constants increase quite steadily as functions of $p$. The plot is that of a convex function, albeit one that is quite close to being linear. Note that $E(\gamma)$ is also increasing in $p$, so this matches the observation from Equation~\eqref{eq:const_inc_E} at the end of the previous section that the constants are increasing in $E(\gamma)$. Secondly, in the interval $[0.9, 1]$ the increasing trend sharpens quite considerably. This suggests that the case $p \rightarrow 1$ is worth looking into in more detail.

\subsection{p tends to 1}\label{subsec:p_to_1}

In this section, we analyse in detail the limit case where $p$ tends to $1$. It is worth noting that here we are dealing with two simultaneous limits, which are the number of grains $N$ (tending to infinity), and the probability of the Binomial distribution $p$ (tending to $1$). As such, we will consider $p = p_N$ to be a function of $N$ such that $p_N \rightarrow 1$ as $N \rightarrow \infty$.

First, to convince ourselves that this case really is worth studying, we ``zoom in'' on the behaviours of the radius and avalanche constants in the interval $[0.9, 1]$. As in Section~\ref{subsec:general_p}, we take $M=100$ and $N=10^7$ for our numerical estimates. The plot is shown on Figure~\ref{fig:radius_avalanche_bin_0.9to1}.

\begin{figure}[h]
\centering
\includegraphics[width=0.9\textwidth]{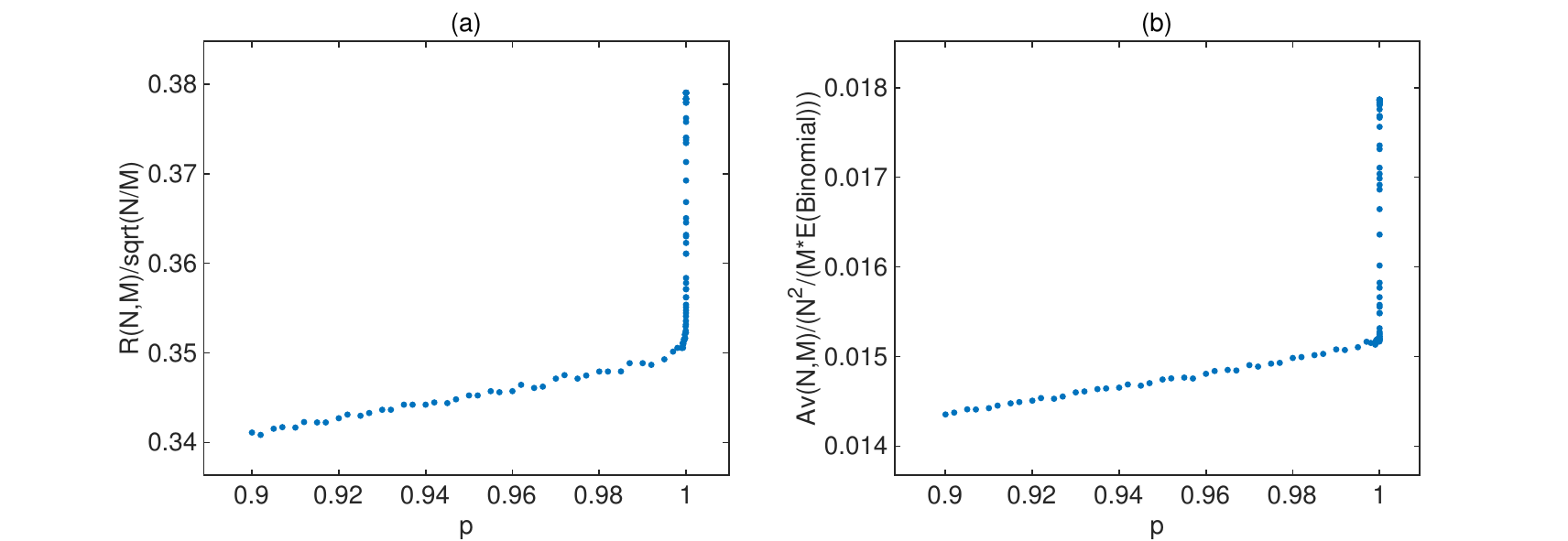}
\caption{Estimates for the radius constant (a) and avalanche constant (b) from $p=0.9$ to $p=1$ when $M=100$ and $N=10^7$}
\label{fig:radius_avalanche_bin_0.9to1}
\end{figure}

These plots clearly indicate that there is a ``jump'' in the values of the constants as $p$ gets close to $1$. We can think of this as a phase transition between the standard binomial regime and the Deterministic model (which corresponds to $p = 1$). Furthermore, experimental data suggests that this phase transition occurs at the scale $1-p \approx \frac{1}{N}$. Note that a phase transition at the same scale of $\frac{1}{N}$ is typical of a number of other probabilistic models, such as the Erd\"os-Renyi random graphs (see e.g.~\cite{JLR2000}). We now provide some heuristic arguments for why the phase transition occurs at this scale in our Binomial SSM.

For clarity of notation, we write $p_N = p$ to indicate that $p$ depends on $N$, and set $q_{N}:=1-p_{N}$. Throughout, we use $C$ to represent a positive finite constant, independent of $N$ and $M$, which can take different values in different equations.

To start, consider a single toppling at a given vertex. The probability that this is a ``full'' toppling, i.e. that $\gamma = M$, is just
\begin{equation}
P(\FT) := P(\gamma=M)=P(\mathrm{Bin}(M,p_{N})=M)=\dfrac{{p_{N}}^M}{1 - {q_N}^M}. \label{eq:probability_full_toppling}
\end{equation}
Here, writing $p_{N}=1-q_{N}$, and taking standard estimates gives
\begin{equation}
P(\FT) \approx 1 - M \cdot q_N. \label{eq:probability_full_toppling_std_est}
\end{equation}
Therefore, the probability of a non-full toppling occurring is
\begin{equation}
P(\NFT) = 1-P(\FT) \approx M\cdot q_{N}. \label{eq:probability_nonfull_toppling}
\end{equation}

Switch the observation perspective into how many non-full topplings there are in the stabilisation of the initial configuration of $N$ grains dropped at the central vertex. From Equation~\eqref{eq:def_avalanche}, we infer that the total number of topplings is roughly $C \cdot \frac{N^2}{(M \cdot E(\gamma))}$ (this is the avalanche number). Thus, the expected value for the number of non-full topplings can be estimated by:
\begin{equation}
\begin{split}
E(\text{total } \NFT \text{ in stabilisation})& = E \big( Av(\mathrm{Bin}; N, M) \big) \cdot  P(NFT) \\
&  \approx C \cdot \frac{N^2}{(M \cdot E(\gamma))} \cdot \left( M \cdot q_N \right) \\
& \approx C \cdot \frac{N^2}{E(\gamma)}\cdot q_{N}, \label{eq:E_total_nonfull_toppling}
\end{split}
\end{equation}
using the estimate from Equation~\eqref{eq:probability_nonfull_toppling} above.

Finally, recall from Equation~\eqref{eq:def_radius} that the radius number is of order $\sqrt{\frac{N}{M}}$. This means that there are roughly $C\cdot \frac{N}{M}$ vertices which are reached in the stabilisation process. Moreover, in this stabilisation, all the reached vertices will topple, except perhaps some on the boundary of the ``toppling ball''. Such non-toppling vertices number at most $\sqrt{\frac{N}{M}} \ll \frac{N}{M}$. Therefore, there are roughly $C\cdot \frac{N}{M}$ vertices which topple in the stabilisation process. As such, if taking a typical vertex $v$ in the toppling ball, its average number of non-full topplings can be estimated by
\begin{equation}
\begin{split}
E(\NFT\text{ at } v  \text{ in stabilisation})\approx \frac{E(\text{total } \NFT \text{ in stabilisation})}{C\cdot (N/M)}\\
\approx C\cdot \frac{\left(\frac{N^2}{E(\gamma)}\right)}{\left(\frac{N}{M}\right)} \cdot q_{N} \approx C \cdot \frac{M}{E(\gamma)} \cdot \left(Nq_{N}\right), \label{eq:E_nonfull_toppling_at_v}
\end{split}
\end{equation}
using estimates from Equation~\eqref{eq:E_total_nonfull_toppling} above.

Now recall that in the (conditioned) binomial case, $E(\gamma)=\frac{p_{N} \cdot M}{(1-(1-p_{N})^M)}$, so that:
\begin{equation} \label{eq:M/E(Binomial)}
\frac{M}{E(\gamma)}=\frac{\left(1-(1-p_{N})^M\right)}{p_{N}}=\frac{\left(1-{q_{N}}^M\right)}{(1-q_{N})}.
\end{equation}
This can be simplified to
\begin{equation} \label{eq:simplication_M/E(Binomial)}
\frac{M}{E(\gamma)}=1+q_{N}+{q_{N}}^2+\cdots+{q_{N}}^{(M-1)}.
\end{equation}
As $N$ tends to infinity, $q_{N}$ tends to 0, so that we get:
\begin{equation} \label{eq:small-o_M/E(Binomial)}
\frac{M}{E(\gamma)}=1+o(1).
\end{equation}
This means that ultimately the multiplication by $\frac{M}{E(\gamma)}$ in Equation~\eqref{eq:E_nonfull_toppling_at_v} can be ignored, and therefore the equation becomes:
\begin{equation} \label{eq:final_E_nonfull_toppling_at_v}
E(\NFT \text{ at } v  \text{ in stabilisation}) \approx C \cdot \left( Nq_{N} \right).
\end{equation}

Considering Equation~\eqref{eq:final_E_nonfull_toppling_at_v} and our empirical data, what seems to happen is that the limit behaves differently depending only on the product $N\cdot q_{N}$. That is, we distinguish three cases, depending on whether the average number of non-full topplings at a vertex $v$ (where $v$ is a typical vertex of the toppling ball) is 0, finite, or infinite.

\begin{itemize}

\item \textbf{Case~1:} $Nq_N \rightarrow 0$, i.e. $q_{N}\ll \frac{1}{N}$. In this case, at a typical vertex $v$, there are only full topplings on average, and the model behaves like the Deterministic one. 

\item \textbf{Case~2:} $Nq_N \rightarrow + \infty$, i.e. $q_{N}\gg \frac{1}{N}$. In this case, at a typical vertex $v$, there are infinitely many non-full topplings on average, and the model behaves like a “standard binomial” one.

\item \textbf{Case~3:} $Nq_N \rightarrow a$, i.e. $q_{N}\sim \frac{a}{N}$, for some finite $a \in (0, +\infty)$. In this case, a typical vertex $v$ topples non-fully a positive, but finite, number of times on average, and the model behaves neither like the Deterministic model, nor like a standard binomial one, but like something in between.

\end{itemize}

This indicates that we have three possible regimes to consider: $1-p \ll \frac{1}{N}$, $1-p \gg \frac{1}{N}$, and $1-p=\frac{a}{N}$.

\subsubsection{The regime $1-p \ll \frac{1}{N}$}\label{subsubsec:p_first_then_N}

There are two interpretations possible for the $1-p\ll \frac{1}{N}$ regime. The standard interpretation is that “$p$ tends to $1$ much faster than $N$ tends to infinity”, but it could also be thought of as “$p$ tends to $1$ first, then $N$ tends to infinity”.

We provide data for two examples of families which belong to this regime: $1-p=\left(\frac{1}{N}\right)^{\frac{3}{2}}$ and $1-p=\left(\frac{1}{N}\right)^2$,  and provide numerical estimates for the radius and avalanche constants in both examples. Figures~\ref{fig:radius_avalanche_bin_3^2} and \ref{fig:radius_avalanche_bin_2} plot the numerical estimates of the radius and avalanche constants for these two examples. The same quotients as in Section~\ref{subsec:numeric_estimates} are plotted as functions of $N$ to observe the limiting behaviour. In both figures, we consider the values $M=10$ and $M=100$, with $N$ ranging from $10^4$ to $10^6$. It is important to remember that the parameter $p$ is a function of $N$ here, so that different data points in $N$ also correspond to different values of $p$.

\begin{figure}[h]
\centering
\includegraphics[width=0.9\textwidth]{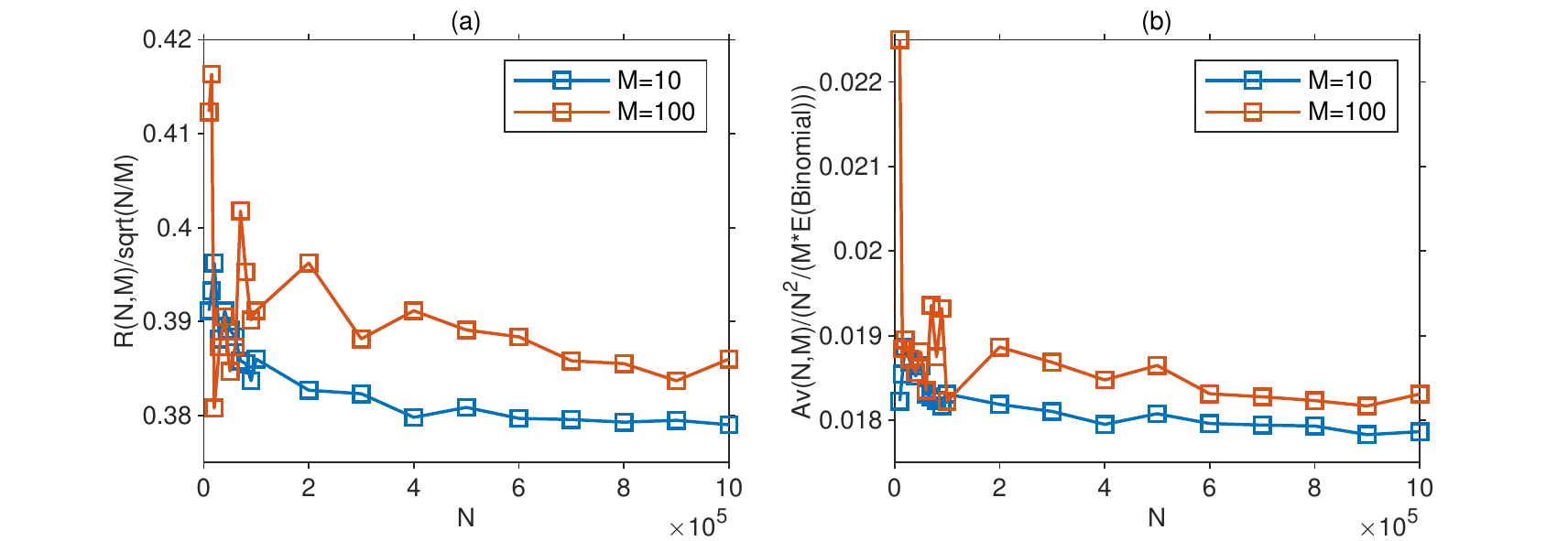}
\caption{Estimates for the radius constant (a) and avalanche constant (b) for $1-p=\left(\frac{1}{N}\right)^{\frac{3}{2}}$ for $M=10$ and $M=100$}
\label{fig:radius_avalanche_bin_3^2}
\end{figure}

\begin{figure}[h]
\centering
\includegraphics[width=0.9\textwidth]{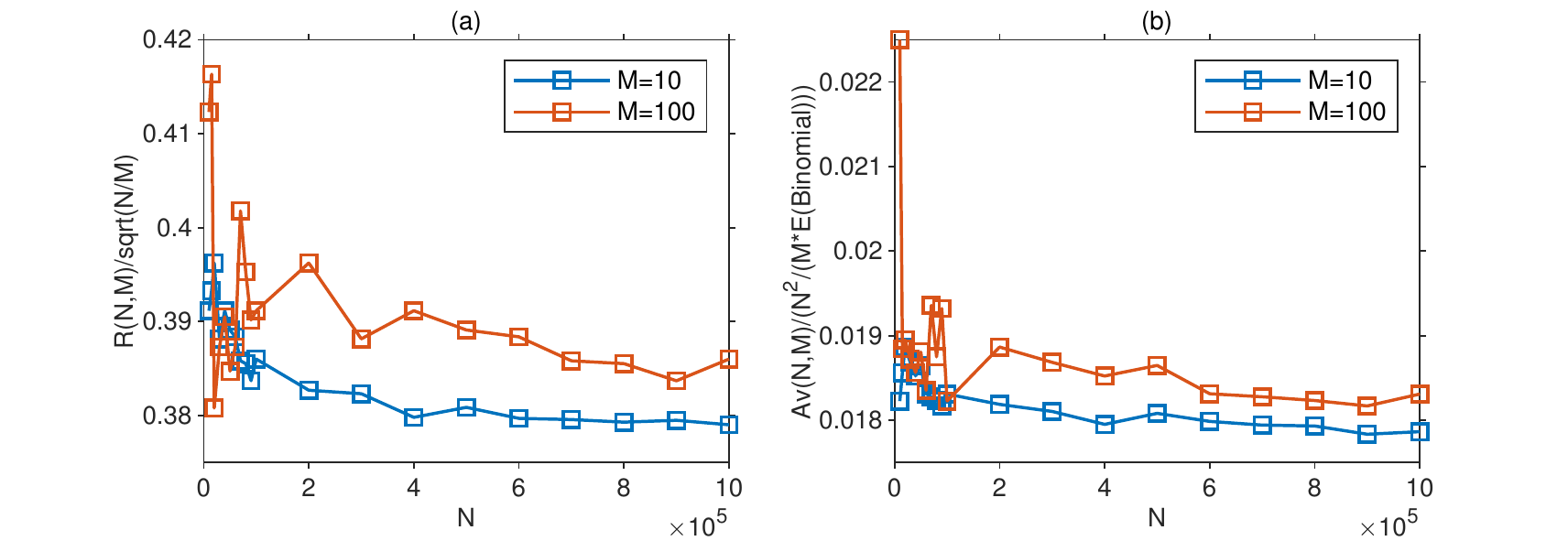}
\caption{Estimates for the radius constant (a) and avalanche constant (b) for $1-p=\left(\frac{1}{N}\right)^2$ for $M=10$ and $M=100$}
\label{fig:radius_avalanche_bin_2}
\end{figure}

We observe a striking similarity between the two figures, i.e. the cases $1-p=\left(\frac{1}{N}\right)^{\frac{3}{2}}$ and $1-p=\left(\frac{1}{N}\right)^2$, for both radius and avalanche constants. These two plots are also remarkably similar to the data for the Deterministic model exhibited in Figure~\ref{fig:radius_avalanche_Det}. This confirms our interpretation that when $1-p \ll \frac{1}{N}$, we are essentially in the Deterministic regime. In particular, the radius constant is roughly $0.38\cdots$, and the avalanche constant $0.018\cdots$.

\subsubsection{The regime $1-p \gg \frac{1}{N}$}\label{subsubsec:N_first_then_p}

Similarly to the $1-p \ll \frac{1}{N}$ regime, two different interpretations can be given to $1-p \gg \frac{1}{N}$: it could be regarded as “$p$ tends to $1$ much slower than $N$ tends to infinity”, or it could be thought of as “$N$ tends to infinity first, then $p$ tends to $1$”. 

To illustrate this regime, we provide data for $1-p=\sqrt{\frac{1}{N}}$, $1-p=\left(\frac{1}{N}\right)^{\frac{1}{3}}$ and $1-p=\left(\frac{1}{N}\right)^\frac{2}{3}$. Figures~\ref{fig:radius_avalanche_bin_1/2}, \ref{fig:radius_avalanche_bin_1/3} and \ref{fig:radius_avalanche_bin_2/3} show the numerical estimates of the radius and avalanche constants, again for $M=10$ and $M=100$, and $N$ ranging from $10^4$ to $10^6$.

\begin{figure}[h]
\centering
\includegraphics[width=0.9\textwidth]{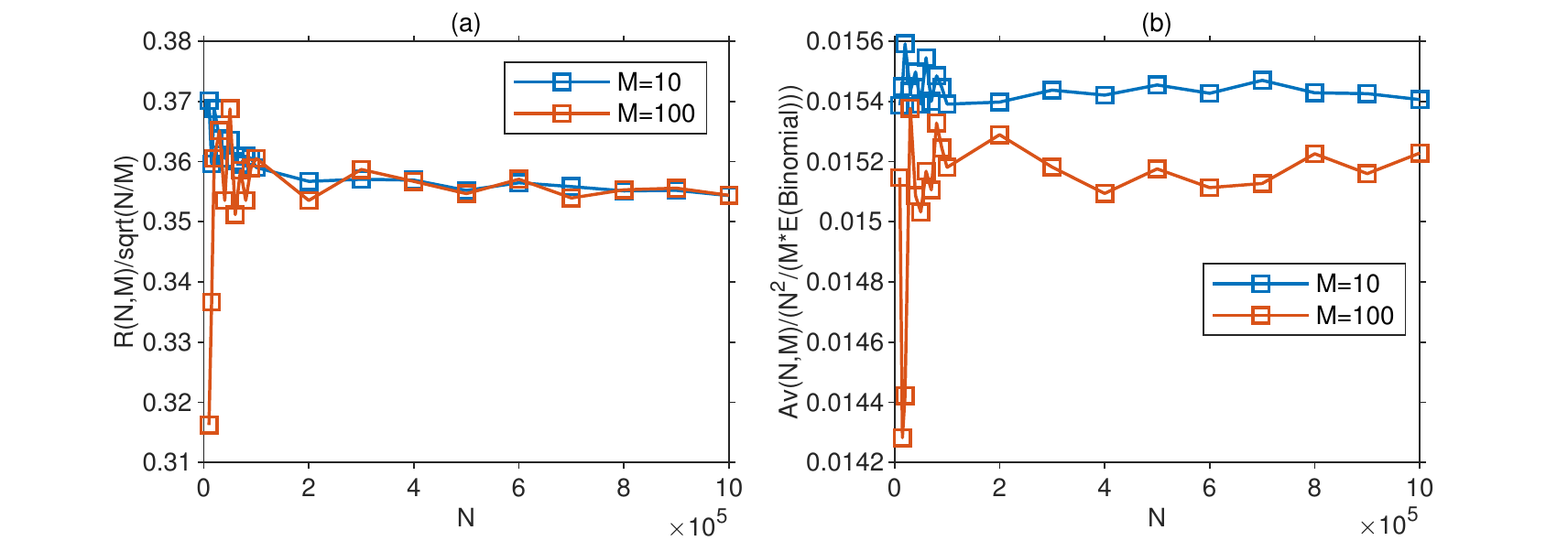}
\caption{Estimates for the radius constant (a) and avalanche constant (b) for $1-p=\left(\frac{1}{N}\right)^{\frac{1}{2}}$ for $M=10$ and $M=100$}
\label{fig:radius_avalanche_bin_1/2}
\end{figure}

\begin{figure}[h]
\centering
\includegraphics[width=0.9\textwidth]{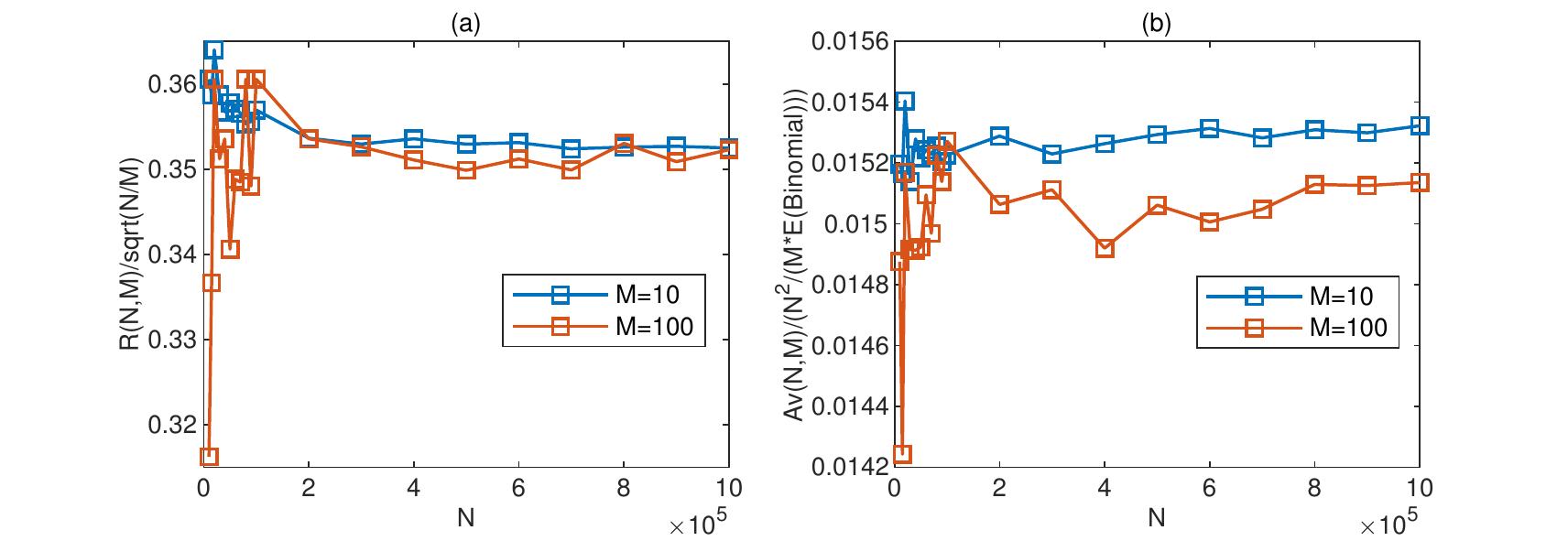}
\caption{Estimates for the radius constant (a) and avalanche constant (b) for $1-p=\left(\frac{1}{N}\right)^{\frac{1}{3}}$ for $M=10$ and $M=100$}
\label{fig:radius_avalanche_bin_1/3}
\end{figure}

\begin{figure}[h]
\centering
\includegraphics[width=0.9\textwidth]{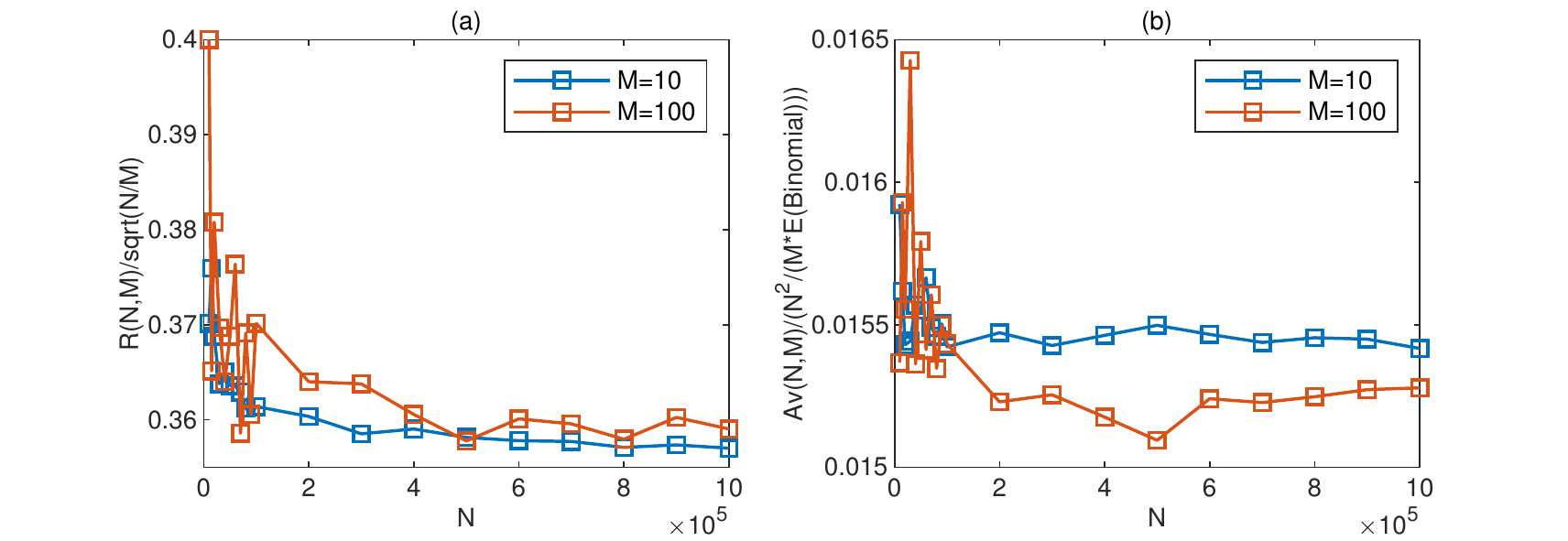}
\caption{Estimates for the radius constant (a) and avalanche constant (b) for $1-p=\left(\frac{1}{N}\right)^{\frac{2}{3}}$ for $M=10$ and $M=100$}
\label{fig:radius_avalanche_bin_2/3}
\end{figure}

The simulations suggest that there is indeed a limit for the radius  and avalanche constants in this regime $1-p \gg \frac{1}{N}$ which does not depend on $M$. Moreover, the numeric value of the radius constant is roughly $0.35\cdots$, while that of the avalanche constant is roughly $0.015\cdots$. These values appear to match the limiting values observed on Figure~\ref{fig:radius_avalanche_bin_0.9to1} before the “jump” as $p$ gets close to $1$. This supports our interpretation of a ``standard binomial'' model for this regime.

\subsubsection{The regime $1-p = \frac{a}{N}$}\label{subsubsec:phase_transition}

In this section, we focus on the regime lying in between that of Sections~\ref{subsubsec:p_first_then_N} and \ref{subsubsec:N_first_then_p}. For this, we introduce an additional parameter $a \in (0, +\infty)$ and set $1-p = \frac{a}{N}$.

Figure~\ref{fig:radius_avalanche_bin_1} gives numerical estimates for the radius and avalanche constants for $a=1$ (i.e. $p = 1 - \frac{1}{N}$), with $N$ ranging from $10^4$ to $10^7$, for both $M=10$ and $M=100$. We observe that the ratios do indeed converge to some limit value, which appears to be independent of $M$. The numeric estimates obtained for the radius and avalanche constants are roughly $0.37$ and $0.017$ respectively, which are in between the values of the $1-p \gg \frac{1}{N}$ and $1-p \ll \frac{1}{N}$ regimes from the previous parts.

\begin{figure}[h]
\centering
\includegraphics[width=0.9\textwidth]{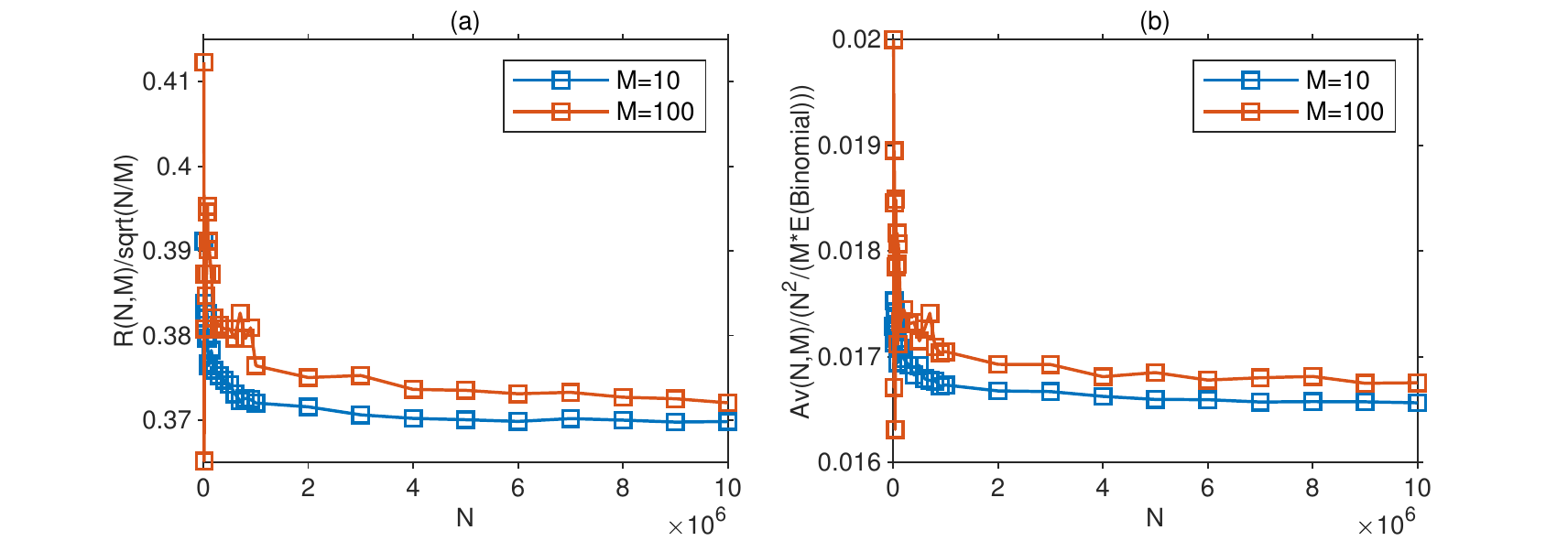}
\caption{Estimates for the radius constant (a) and avalanche constant (b) for $1-p=\frac{1}{N}$ for $M=10$ and $M=100$}
\label{fig:radius_avalanche_bin_1}
\end{figure}

We are also interested in how these constants behave in $a$ when $1-p = \frac{a}{N}$. Figure~\ref{fig:radius_avalanche_bin_a/N} therefore plots the numerical estimates of these constants as functions of $a$ for $a$ ranging from $10^{-5}$ to $10^5$. We consider the values $M=10$ and $M=100$, and take $N = 10^7$. For convenience, we use a log-scale for our X-axis.

\begin{figure}[h]
\centering
\includegraphics[width=0.9\textwidth]{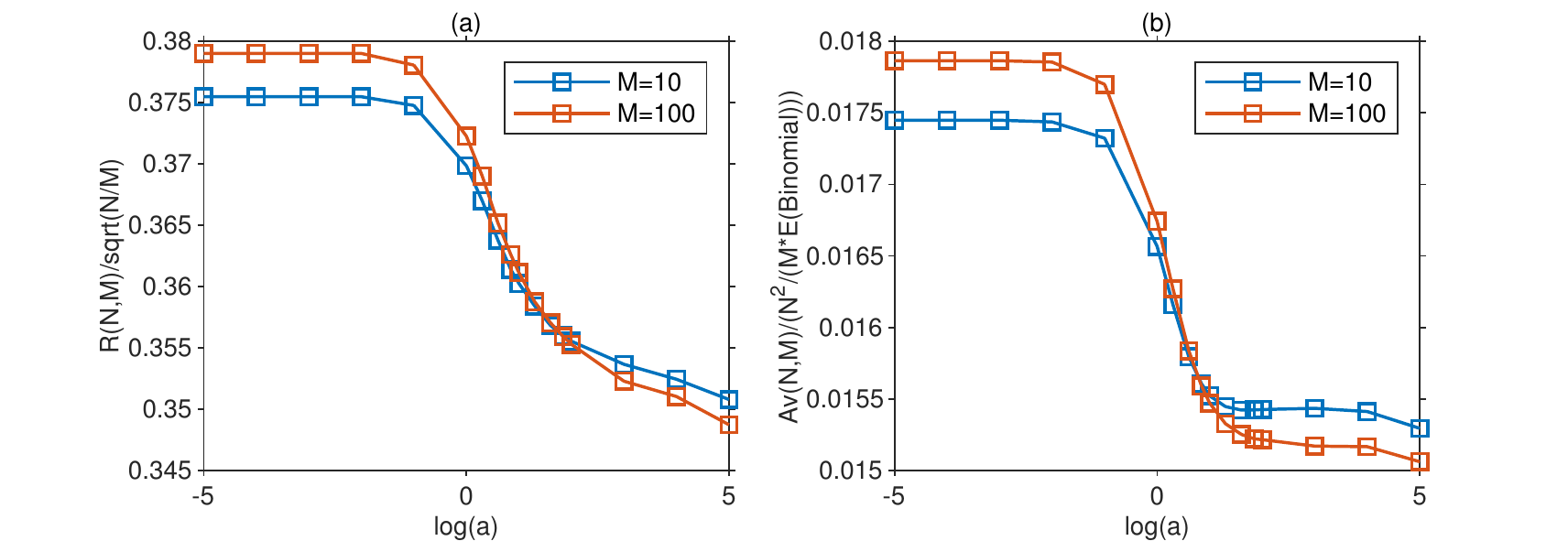}
\caption{Estimates for the radius constant (a) and avalanche constant (b) for $1-p=\frac{a}{N}$ when $\log{(a)}=[-5,5]$ for $M=10$ and $M=100$}
\label{fig:radius_avalanche_bin_a/N}
\end{figure}

Based on our observations, as expected, it seems that the constants depend only on $a$ (and not on $M$). Moreover, as $a$ tends to $0$, the radius constant tends to roughly $0.38$, while the avalanche constant tends to approximately $0.018$, both matching the corresponding values of the Deterministic distribution. Similarly, when $a$ tends to $+ \infty$, the radius and avalanche constants tend to $0.35$ and $0.015$ respectively, coinciding with the values obtained in the regime $1-p \gg \frac{1}{N}$.  Consequently, from the above statements, our interpretation of the regime $1-p = \frac{a}{N}$ which acts as a phase transition between the ``standard binomial'' and Deterministic limits described in  Sections~\ref{subsubsec:p_first_then_N} and \ref{subsubsec:N_first_then_p} is confirmed.

\subsection{p tends to 0}\label{subsec:p_to_0}

We first note that in our model, the probability of $p=0$ is not a possible value (that is, nothing would ever happen). However, if we consider a single toppling, when $p$ tends to 0, the conditioned Binomial distribution is nearly certain to take the value $1$, i.e., one grain gets sent to each neighbour. If this happens every time, then the model simply corresponds to the toppling rules in the Always-1 distribution. This suggests that we may see the Always-1 distribution appear as the limiting case when $p \rightarrow 0$.

Therefore, in Figure~\ref{fig:radius_avalanche_bin_0.1to0}, we use the Always-1 values as the data point for $p=0$. In addition, we take $p$ ranging from $0.006$ to $0.1$. The lower value $p=0.006$ is due to time restrictions, as our simulations run very slowly for low values of $p$ (our code could doubtless be improved here). We then plot estimates of the radius and avalanche constants as functions of $p$ for $N=10^7$ and the values $M=10$, $100$ and $1000$ respectively.

\begin{figure}[h]
\centering
\includegraphics[width=0.9\textwidth]{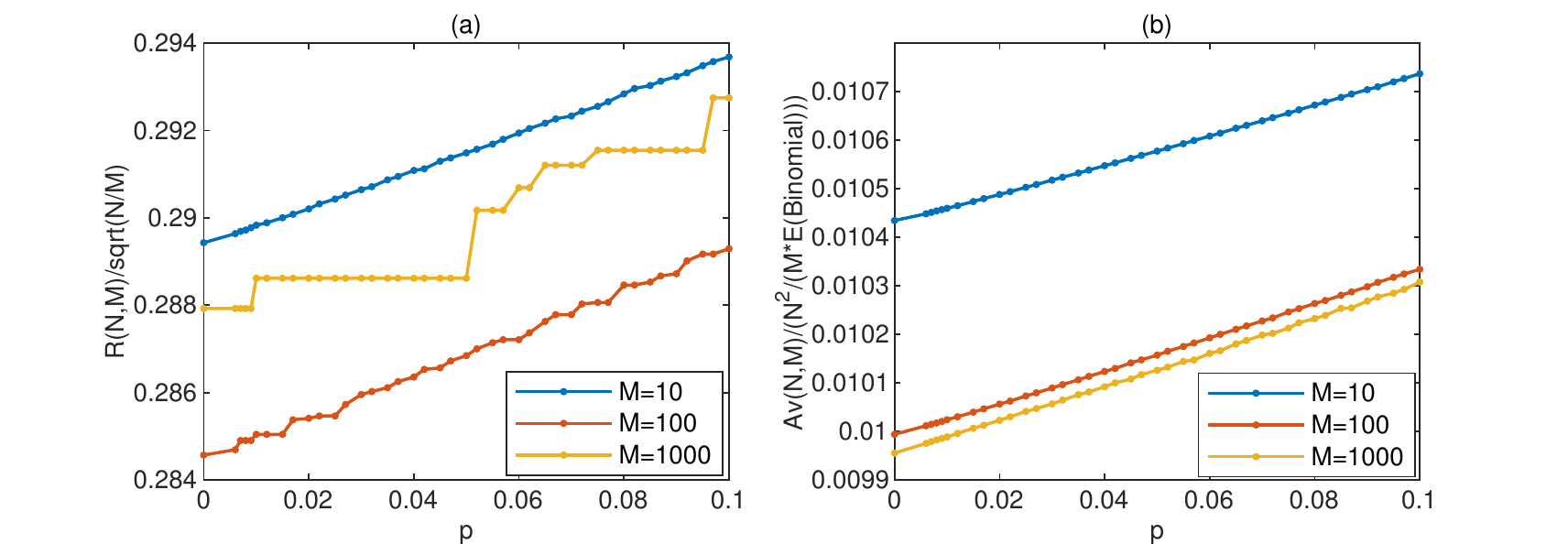}
\caption{Estimates for the radius constant (a) and avalanche constant (b) from $p=0.1$ to $p=0$ for $M=10$ and $M=100$}
\label{fig:radius_avalanche_bin_0.1to0}
\end{figure}

The experimental data matches our interpretation above of the Always-1 distribution as the limiting case of the Binomial model as $p$ tends to $0$. The radius and avalanche constants both converge quite smoothly to the estimates of the constants for the Always-1 distribution. The radius constant in the case $M=1000$ does behave a little strangely, but this is probably because this is a large value of $M$, roughly comparable to $\frac{1}{p}$ in our data: in this case we know that a binomial random variable with parameters $M, p$ behaves like a Poisson random variable, which may explain the non-smoothness of this particular curve.

\section{Conclusion and future works}\label{sec:conclusion_future}

\subsection{Conclusion}\label{subsec:conclusion}

In conclusion, this paper presents an empirical study of the single-source sandpile model for a slight variant of the SSM introduced by Kim and Wang in~\cite{KW2020}. In this model, vertices topple when they have at least $4M$ grains of sand, for some fixed parameter $M$. When toppling, a vertex chooses a random number $k \in \{1, \cdots, M\}$ according to some probability distribution $\gamma$, and sends $k$ grains to each of its four neighbours (in $\Z^2$). The single-source sandpile model consists in initially placing $N$ grains of sand at the origin of $\Z^2$, and making these (random) topplings until no more vertices can topple. We then study the asymptotics of this stabilisation process as $N$ tends to infinity.

In this paper, we have focused on two global parameters of the system. The first is the radius number $R(\gamma; N, M)$ which measures the maximum distance from the origin to which grains are sent in the stabilisation. The second is the avalanche number $Av(\gamma; N, M)$ which counts the total number of topplings during the stabilisation.

We study these two numbers for six different distributions $\gamma$: Deterministic, Always-1, Uniform, Binomial, Log-Law, and Power-Law. Our simulations show the following asymptotics as $N$ tends to infinity:
\begin{equation}\label{eq:asympt_radius_conc}
R(\gamma; N, M) \sim c_R(\gamma) \cdot \sqrt{\frac{N}{M}},
\end{equation}
\begin{equation}\label{eq:asympt_avalanche_conc}
Av(\gamma; N, M) \sim c_{Av}(\gamma) \cdot \dfrac{N^2}{M \cdot E(\gamma)},
\end{equation}
where $c_R(\gamma)$ and $c_{Av}(\gamma)$ are constants depending only on the toppling distribution $\gamma$ and neither on $N$ nor on $M$. We call these the radius and avalanche constants of the system, and give numerical estimates of their values for the six $\gamma$ distributions. We also give heuristic arguments for the behaviours observed in Equations~\eqref{eq:asympt_radius_conc} and \eqref{eq:asympt_avalanche_conc}.

Finally, we analyse in more detail the behaviour of these two constants in the case where $\gamma$ is a Binomial distribution with parameters $M, p$, viewing them as functions of $p \in (0, 1]$. We see that the two constants are increasing in $p$, and that as $p$ tends to $0$ they converge to the constants of the Always-1 distribution, as expected. In the $p \rightarrow 1$ case, we exhibit a phase transition at the scale $1-p \propto \frac{1}{N}$. More precisely, we identify three regimes.
\begin{itemize}
\item If $1-p \ll \frac{1}{N}$, the model behaves like the Deterministic distribution.
\item If $1-p \gg \frac{1}{N}$, the model behaves like a ``standard binomial'' distribution.
\item If $1-p \sim \frac{a}{N}$ for $a > 0$, the model behaves differently. In this case, we give numeric estimates for the radius and avalanche constants for various values of $a$ in the range $[10^{-5}, 10^5]$, and plot these as a function of $a$. The limits as $a$ tends to $0$ and $+ \infty$ give the above two regimes respectively.
\end{itemize}
We give heuristics for why this is the correct scale for this phase transition to occur, based on estimates for the expected number of non-full topplings (where $\gamma < M$) at a typical vertex.

\subsection{Future Work}\label{subsec:future}

We conclude by giving a number of possible future work directions on this topic. Since this is, to our knowledge, the first study of a single-source SSM, there remains much work to be done.

Firstly, we should seek to make our heuristic arguments from Sections~\ref{subsec:heuristics_asymptotics} and \ref{subsec:p_to_1} more rigourous. In particular, it is not completely clear to us why the correct scale for the phase transition when $p$ tends to $1$ depends on the average number of non-full topplings at a given vertex.

We could also analyse the radius and avalanche constants in the Power-Law case as functions of the parameter $s>0$. This would require us to improve our code (simulations for this distribution were very slow). In particular, we would expect the model to behave like the Uniform model when $s$ tends to $0$, and the Always-1 model when $s$ tends to infinity (as $s$ tends to infinity, we have $P(\mathrm{Power-Law} = k) \ll  P(\mathrm{Power-Law} = 1)$ for $k \geq 2$, the same behaviour as the Binomial distribution when $p$ tends to $0$). It would be interesting to see if experimental data confirms this, and whether there are any phase transitions in these limits.

Another possible avenue of future research would be to focus on more local parameters, rather than the (global) radius and avalanche numbers. For example, we could look at the odometer function, as has been done in the single-source ASM~\cite{PS2013}. We could also try to describe the limit ``shapes'' of the stable configuration reached after all topplings are complete. Figure~\ref{fig:shapes} gives some examples of these, and it would be interesting to analyse them further.

At first glance, the shapes in Figure~\ref{fig:shapes} show some interesting properties. Recall that we use a sliding RGB colour scale, with blue~$\equiv 1$ grain, green~$\equiv 2M$, and red~$\equiv 4M-1$. We notice that they are largely made up of a mix of red and green. This indicates that the number of grains at any vertex in the toppling ball would be at least $2M$. We also see that the proportion of ``red'' vertices, i.e. vertices with closer to $4M-1$ grains (vertices that are ``nearly full''), is decreasing in $E(\gamma)$. For example, the Always-1 distribution has almost all full vertices, whereas the $\mathrm{Binomial}(p = 0.75)$ distribution, which has the highest expected value among the six examples of the figure, has the largest proportion of green (half-full) vertices. This matches with findings from Equation~\eqref{eq:const_inc_E} that the radius increases with $E(\gamma)$. Indeed, the balls all have $N$ grains in total, since the total number of grains in the system remains constant through the stabilisation process. This means that if the ball is larger (has higher radius), we should expect individual vertices to have fewer grains. It would be interesting to formalise these observations, and also to study the shapes more generally.

Finally, we could attempt to extend our work to study the single-source SSM with one or more sinks. A  sink is introduced to allow sand grains to exit the system, and such a model was studied for the single-source ASM in~\cite{SD2010}. In particular, the authors showed that sinks will change the radius scale when initially placing $N$ sand grains at the origin: if there is a sink, the radius behaves like $\sqrt{N/\log{N}}$ (for large $N$); while, as we have seen, if there is no sink, the radius will behave like $\sqrt{N}$. It would be interesting if a similar change of scale occurs for the SSM.

\bibliographystyle{plain}
\bibliography{SingleSource}

\end{document}